\documentclass[a4paper,12pt]{article}

\usepackage{amsthm,amsmath,amscd,mathtools}
\usepackage[vmargin=2cm,hmargin=2cm,headheight=14.5pt,top=2cm,headsep=.5cm]{geometry}
\usepackage{hyperref}
\usepackage[utf8]{inputenc}
\usepackage[autobold]{mathfixs}
\usepackage{stackrel}
\usepackage{cases}
\usepackage{circuitikz}
\usepackage{tikz}
\usepackage{placeins}

\usepackage[charter]{mathdesign}
\usepackage{bm}
\usepackage{caption}
\usepackage{subcaption}
\usepackage{xcolor}
\usepackage{graphicx}
\DeclareMathSizes{12}{12}{8}{6}

\usepackage{cite}

\usepackage[shortcuts]{extdash}
\usepackage{xurl}
\usepackage{enumitem}
\setlist[itemize,2]{label=$\centerdot$}
\setlist[itemize,3]{label=$\triangle$}

\newtheoremstyle{ptheorem}{1em}{0em}{\itshape}{}{\bfseries}{.}{.5em}{\thmname{#1}\thmnumber{
		#2}\thmnote{ (\hspace{-.01pt}{#3})}}

\theoremstyle{ptheorem}

\newtheorem{thm}{Theorem}[section]
\newtheorem{pro}[thm]{Proposition}
\newtheorem{lem}[thm]{Lemma}
\newtheorem{cor}[thm]{Corollary}

\newtheoremstyle{hdef}{1em}{0em}{}{}{\bfseries}{.}{.5em}{\thmname{#1}\thmnumber{
		#2}\thmnote{ (\hspace{-.01pt}{#3})}}
\theoremstyle{hdef}

\newtheorem{dfn}[thm]{Definition}
\newtheorem{rem}[thm]{Remark}

\numberwithin{equation}{section}
\numberwithin{figure}{section}

\DeclareMathOperator{\dif}{d}

\newcommand{\bR}{{\mathbb R}}

\renewcommand{\phi}{\varphi}
\renewcommand{\le}{\leqslant}
\renewcommand{\ge}{\geqslant}

\newcommand{\bs}{\backslash}

\newcommand{\olb}[1]{%
	\vbox{\offinterlineskip\ialign{\hfil##\hfil\cr
			$\rotatebox[origin=c]{90}{$]$}$\cr\noalign{\kern-.45ex}{$#1$}\cr}}}

\newcommand{\noop}[1]{}

\usepackage{stmaryrd}

\def\arraystretch{1.2}
\allowdisplaybreaks

\begin{document}

\title{Stieltjes differential equations with bounded-variation derivators and application to thermal stress in solar panels}
\author{Lamiae Maia$^1$ and F. Adri\'an F. Tojo$^2$}
\date{}
	\maketitle

\begin{center} {\small $^{1}$ LAMA Laboratory, Mathematics Department, Faculty of Sciences\\
		Mohammed V University in Rabat,  Rabat,  Morocco. \\ e-mail: \emph{lamiae\_maia@um5.ac.ma}}
	\\
	\small $^{2}$ Departamento de Estat\'{\i}stica, An\'alise Matem\'atica e Optimizaci\'on \\ Universidade de Santiago de Compostela \\ 15782, Facultade de Matem\'aticas, Campus Vida, Santiago, Spain.\\ CITMAga, Santiago de Compostela \\ e-mail: \emph{fernandoadrian.fernandez@usc.es}
\end{center}

	\medbreak

	\begin{abstract}
	In this work we extend the theory of Stieltjes systems beyond the monotone case, establishing new chain rules, generalized versions of the Fundamental Theorem of Calculus, compactness tools for Peano–type results, and a $g$–exponential for explicit linear solutions. Continuity notions relative to vector-valued derivators further allow us to study everywhere differentiable solutions. As an application, we model thermal stress effects on solar panels and battery health, highlighting the practical value of non-monotone derivators.
	\end{abstract}
 \maketitle

\medbreak


\medbreak


\section{Introduction}
In recent years, Stieltjes differential equations have attracted growing attention as a flexible framework for modeling dynamical systems presenting jumps or stationary periods (see, e.g.~\cite{AFNT,FernandezTojoStieltjesBochnerSpaces,PM,PM2,PM3,MEF1,MEF2,MEF3,MaiaTojo-1-2025}) that classical differential models may struggle to describe adequately. These equations generalize the classical theory by replacing the standard derivative with one taken with respect to a function of bounded variation—the so-called {\it derivator}.  Most of the existing literature has focused on the case of left-continuous monotone derivators~\cite{FernándezAlbésTojo2024FirstNSec,FP,PR}. In contrast, the study of Stieltjes differential equations driven by derivators with weaker regularity remains relatively scarce. An important step in this direction was taken in~\cite{FrTo}, the monotonicity condition has been relaxed by considering derivators with {\it controlled variation} within the class of locally bounded variation functions. These functions display monotone behavior on each connected component, which is determined by excluding a negligible set. The exploration of Stieltjes differential equations beyond the monotone framework is both mathematically intriguing and practically useful, since it provides a natural bridge between Stieltjes and measure differential equations~\cite{MonteiroBianca-GenerDeriv2017}.

More recently, in \cite{MaiaTojo-1-2025}, the authors carried the analysis beyond the monotone case, establishing results under no monotonicity assumptions. These results include new generalized versions of the Fundamental Theorem of Calculus for Lebesgue–Stieltjes integrals presenting both "almost-everywhere" and "everywhere" versions. The latter is shown to be optimal by means of a new density result~\cite[Theorem 5.10]{MaiaTojo-1-2025}. Building on these results, the present work investigates systems of Stieltjes differential equations of the form
\begin{equation}\label{eq:introduction:SDE}
  \begin{aligned}
  \mathbf{x}_\mathbf{g}'(t)=&\mathbf{f}(t,\mathbf{x}(t)),\quad\text{ for $\mathbf{g}$-almost every }t\in [t_0,t_0+T],\\
  \mathbf{x}(0)=&\mathbf{x}_{0}\in \mathbb{R}^n,
  \end{aligned}
\end{equation}
where $\mathbf{f}=(f_1,\dots,f_n):[t_0,t_0+T]\times\mathbb{R}^n\to\mathbb{R}^n$, $\mathbf{g}=(g_1,\dots,g_n):[t_0,t_0+T]\to\mathbb{R}^n$, with $g_i:[t_0,t_0+T]\to\mathbb{R}$ a left-continuous derivator of {bounded} variation. We extend the results in~\cite{MaiaTojo-1-2025}, by providing two formulas for the Stieltjes derivative of a composition of functions in a similar fashion to~\cite{MarSlav2024generChainRule}. Besides, we adapt compactness arguments from~\cite{FP} to the non-monotonic case.  These results will serve as a tool to establish a Peano–type existence result. We then introduce a g-exponential function, generalizing the one introduced in~\cite{FP,FrTo}, which provides explicit solutions to linear problems.
In addition, we formulate an "everywhere" version of the fundamental theorem of calculus involving several non-monotone derivators by means of continuity notions relative to vector-valued derivators. This permits us to study "everywhere" solutions for problem~\eqref{eq:introduction:SDE}, i.e. solutions that are Stieltjes differentiable everywhere, thereby generalizing the monotone case studied in~\cite{FP,MT}. {Moreover}, we extend Picard and Peano type results in the spirit of~\cite{FrTo,PM3}. Finally we apply one of the existence results established earlier to the modeling of thermal stress effects on solar panels performance and battery health, involving derivators of bounded variation.

The paper is structured as follows: In Section~2, we present preliminary tools from Lebesgue\--Stieltjes {integration involving} a left-continuous derivator {of} {bounded} variation. Then we present the Stieltjes derivative with respect to such a derivator and establish two formulations of the Stieltjes derivative of a composition of functions. Section~3 is devoted continuity concepts, including the new versions of the Fundamental Theorem of Calculus involving Lebesgue-Stieltjes integrals in this broader context. In Section~4, we prove compactness results that provide the key tool for establishing a Peano–type existence theorem, and in Section~5 we introduce a $g$–exponential function that permits explicit solutions of linear equations. Section~6 forms the core of this paper. We begin by presenting continuity notions with respect to vector–valued derivators, and then establish an "everywhere" version of the Fundamental Theorem of Calculus involving several derivators.  This will serve as a tool to prove "everywhere" $g$-differentiability of solutions. We then apply the previous results to systems of Stieltjes differential equations with vector–valued derivators. By reducing problem~\eqref{eq:introduction:SDE} to an abstract fixed point problem, we establish Picard and Peano type theorems for problem~\eqref{eq:introduction:SDE} in the Carathéodory sense. Section~7 illustrates the applicability of the theory by modeling the impact of thermal stress on photovoltaic panel performance and the state of health of the battery.

\section{Preliminaries}

Let $[a,b]\subset \mathbb{R}$, and  $g:[a,b]\to {\mathbb R}$ be a function. Consider ${\mathcal P}([a,b])$ the set of the partitions of the interval $[a,b]$ defined by
\[
{\mathcal P}([a,b]):=\{ P=(\tau_1,\dots,\tau_{n_P}):\, n_P\geqslant 2;\, \tau_k\in [a,b],\, k=1,\dots, n_P\, ; \tau_k\leqslant \tau_{k+1};\, k=1,\dots, n_P -1\}.
\]
 The {\em total variation} of $g$ in $[a,b]$ is defined by
\[
\operatorname{var}_g [a,b]=\sup_{P\in {\mathcal P}([a,b])}\sum_{k=1}^{n_P -1}|g(\tau_k)-g(\tau_{k+1})|.
\]
If $\operatorname{var}_g [a,b] < \infty$, then  we say that $g$ is of {\em bounded variation}.

We denote by $\operatorname{BV}([a,b],{\mathbb R})$ the space of real\--valued functions of {\em  bounded variation in $[a,b]$}, and we denote by $\operatorname{\operatorname{BV}^{lc}}([a,b],{\mathbb R})$  the space of functions of $\operatorname{BV}([a,b],{\mathbb R})$  which are also left\--continuous.

Here and afterwards, we fix a derivator $g\in \operatorname{\operatorname{BV}^{lc}}([a,b],{\mathbb R})$. In particular, $g$ is regulated, and we denote{
\[
D_g=\{t\in[a,b):\Delta^+g(t)\neq 0\},
\]}
where $\Delta^+g(t):= g(t^+)-g(t)$ and $g(t^+):=\lim\limits_{s\to t^+} g(s)$. By Froda's Theorem~\cite{Froda}, $D_g$ is at most countable.  Moreover, we define
\begin{equation*}
	C_g:=\{ t \in [a,b] \, : \text{$g(s)=g(t)$ for all $s\in(t-\varepsilon,t+\varepsilon)$ for
		some $\varepsilon>0$} \}.
\end{equation*}
The set $C_g$ is open in the usual topology of $\mathbb R$, thus, it can be written in the form
\begin{equation*}
	C_g=\bigcup_{n \in \Lambda}(a_n,b_n),
\end{equation*}
where
$\Lambda\subset \mathbb{N}$ and $(a_i,b_i)\cap (a_j,b_j)=
\emptyset$ for $i\neq j$. With this notation, we denote
\[
N_g^-:=\{a_n\}_{n \in \Lambda}\backslash (D_g\cup\{a\}),\quad N_g^+:=\{b_n\}_{n \in \Lambda}\backslash D_g,\quad \text{and}\quad
N_g:=N_g^-\cup N_g^+.
\]
{If $b\notin N_g^+$, then for} each $t\in[a,b]$, we consider the following notation introduced in~\cite[Proposition~3.9]{FMarTo-OnFirstandSec}:
\begin{equation}\label{eq:notation of t^*}
 t^*=\begin{dcases}
 b_n, & \mbox{if $t\in(a_n,b_n)\subset C_g$, or ($t=a\notin D_g$ and $(a,b_n) \subset C_g$)},\\
  t, & \mbox{{ot}herwise.}
 \end{dcases}
 \end{equation}
{Observe that $t^* \notin C_g$ for all $t\in [a,b]$.}
\begin{dfn}\label{dfn:variation function}
 We define the {\em variation function} of $g$ by the function {$\widetilde{g}:=\operatorname{var}(g):[a,b]\to{\mathbb R}$ given by
 \[
 \widetilde{g}(t):=\operatorname{var}(g)(t)=\operatorname{var}_g [a,t].
 \]}
\end{dfn}
\begin{rem}\label{rem:var[a,b]=tilde(g)(b)-tilde(g)(a)}
The function $\widetilde{g}$ is nondecreasing, shares the same discontinuity points of~$g$ and its constancy intervals, and for $x,y\in[a,b]$, with $x<y$, we have
\begin{equation*}
 \operatorname{var}_g [x,y]=\widetilde{g}(y)-\widetilde{g}(x).
\end{equation*}
Furthermore, as $g$ is left\--continuous, so is~$\widetilde{g}$.
\end{rem}

Since $g\in  \operatorname{\operatorname{BV}^{lc}}([a,b],{\mathbb R})$, $g$~generates a signed measure, denoted~$\mu_g$, on the measurable { space~$([a,b],\mathcal{M}_g)$, where $\mathcal{M}_g$ is the Borel $\sigma$-algebra induced by the usual topology on $[a,b]$.} The signed measure~$\mu_g$ is defined via the following fundamental formula:
\begin{equation*}
 \mu_g([a,t))=g(t)-g(a),\quad\text{for all }t\in [a,b].
\end{equation*}
The concept of measurability with respect to this signed measure is referred to as $g$\--measurability. {A crucial tool for the analysis of such derivators and their associated measures is the} Jordan decomposition~\cite[Theorem~4.1.2]{BenedettoCzaja2009integration}. This decomposition allows us to express the derivator~$g$ as the difference {of} two monotone, nondecreasing derivators $g_1,g_2\in \operatorname{\operatorname{BV}^{lc}}([a,b],{\mathbb R})$. Each of these derivators is associated with a Lebesgue\--Stieltjes measure. Furthermore, this decomposition {is closely related to } Hahn's decomposition of signed measures~\cite[Theorem~5.1.6]{BenedettoCzaja2009integration}. {In particular, the measure~$\mu_g$ can be expressed in terms of its {\em positive} and {\em negative variations}, denoted by~$\mu_g^+$ and~$\mu_g^{-}$, respectively. Formally} the measure space $([a,b],{\mathcal M}_g,\mu_g)$ can be partitioned into two disjoint subsets~$A_g^+$ and~$A_g^-$ such that:
\begin{enumerate}
 \item $\mu_g(M) \geqslant 0$ for each $g$\--measurable set $M\subset A_g^+$;
 \item $\mu_g(M) \leqslant 0$ for each $g$\--measurable set $M\subset A_g^-$.
\end{enumerate}
With this partition,~$\mu_g^+$ and~$\mu_g^{-}$ are defined as
\[
\mu_g^+(M):=\mu_g(M\cap A_g^+),\text{ and } \mu_g^-(M):=-\mu_g(M\cap A_g^-)\quad \text{for every }M\in {\mathcal M}_g.
\]
We can express the {\em total variation} of the signed measure~$\mu_g$ as $|\mu_g|:=\mu_g^+ + \mu_g^-$, which captures the overall magnitude of~$\mu_g$ regardless of sign.

{If we define the functions}~$g_1$  and~$g_2$ on the interval $[a,b]$ as follows:
\[
g_1(t):=\mu_g^+([a,t))\quad\text{and}\quad g_2(t):=\mu_g^-([a,t)),\quad \text{for every }t\in [a,b],
\]
then~$g_1$  and~$g_2$ are both left\--continuous and nondecreasing, and they satisfy the {relation $g=g(a)+g_1-g_2$.} Moreover, the measures~$\mu_{g_1}$  and~$\mu_{g_2}$, which are the Lebesgue\--Stieltjes measures associated with $g_1$ and~$g_2$, coincide with the positive and negative variations of~$\mu_g$, respectively. It is also noteworthy that the measure~$|\mu_g|$ coincides with the Lebesgue\--Stieltjes measure associated with the variation function~$\widetilde{g}$, the reader is referred to~\cite[Theorems~5.1.6,~5.1.8, and~5.1.9]{BenedettoCzaja2009integration} for further details. This measure satisfies the property:
\[
|\mu_g|([u,v))=\widetilde{g}(v)-\widetilde{g}(u)=\operatorname{var}_g[u,v]\quad\text{for every $u,v\in[a,b]$, with $u<v$.}
\]
{We refer to the \emph{measurability} (resp. \emph{integrability}) with respect to the measure $|\mu_g|$ by $\widetilde{g}$\--measurability (resp. $g$\--integrability), and we denote $L_g^1([a,b),\mathbb{R})$ the space of $g$\--integrable functions on the interval $[a,b) \subset \mathbb{R}$ endowed with the norm
\[
\|f\|_{L_{g}^{1}([a,b))}:=\int_{[a,b)} |f| \operatorname{d} |\mu_g|, \quad \text{ for every } f\in L_{g}^{1}([a,b),\mathbb{R}).
\]
We say that a property holds for \emph{$|\mu_g|$\--almost every $t\in [a,b)$} or  \emph{$|\mu_g|$\--almost everywhere}  (shortly, $|\mu_g|$\--a.e.) if it holds except on a set $N\subset [a,b)$ such that $|\mu_g|(N)=0$.
}

\subsection*{The Stieltjes derivative}
Now, we present the definition of the Stieltjes derivative with respect to non\--monotonic derivators of {bounded} variation, {as  introduced} in~\cite{MaiaTojo-1-2025}. Additionally, we provide two versions of the chain rule, which will be utilized in {Section~5}. Throughout this part, let $[a,b]\subset \mathbb{R}$ be a closed interval and {$g\in \operatorname{\operatorname{BV}^{lc}}([a,b],{\mathbb R})$ such that $b\notin N_g^+$}.
\begin{dfn}\label{dfn:g-derivative(non monotonic case)}
Let $f:[a,b]\to\mathbb R$, $t\in[a,b]$ and assume there exist $d\in{\mathbb R}$, and a function $h:[a,b]\to \mathbb{R}$ satisfying
\begin{enumerate}
 \item $h(t^*)=0$;
 \item $h$ is continuous at $t^*$ if {$t^*\in [a,b]\bs(D_g\cup N_g)$};
 \item $h$ is right\--continuous at $t^*$ if {$t^*\in(D_g\cup N_g^+)$};
 \item $h$ is left\--continuous at $t^*$ if {$t^*\in N_g^-$};
\end{enumerate}
 such that, $f(s)=f(t^*)+[d+h(s)][g(s)-g(t^*)]$ for $s\in[a,b]$, $g(s)\ne g(t^*)$
with $t^*$ as in~\eqref{eq:notation of t^*}. In that case we say that $f'_g(t)\equiv d$ is the {\em Stieltjes derivative} or {\em$g$\--derivative} of $f$ at $t$ and that $f$ is \emph{$g$\--differentiable at $t$}.
\end{dfn}
\begin{rem}\label{rem:ld}{With the above notation, Definition~\ref{dfn:g-derivative(non monotonic case)} can equivalently be written as
	\[ f'_g(t):=\begin{dcases}
		\lim_{\substack{s \to t^*\\ g(s)\ne g(t^*)}}\frac{f(s)-f(t^*)}{g(s)-g(t^*)}, &\mbox{if } t^*\not\in D_{g},\\
		\lim_{s \to t^{*+}}\frac{f(s)-f(t^*)}{g(s)-g(t^*)}, &\mbox{if } t^*\in D_{g},\end{dcases}\]}
cf. \cite[Definition 2.26]{MonteiroBianca-GenerDeriv2017}. {In particular, in the case where $t^*\in D_g$ in Definition~\ref{dfn:g-derivative(non monotonic case)}, since $g$ is regulated, $f'_g(t)$ exists if and only if $f(t^{*+})$ exists,and in that case we have}
\[
f'_g(t)=\frac{f(t^{*+})-f(t^*)}{\Delta^+ g(t^*)}.
\]
\end{rem}

In the following proposition, we state the chain rule for non\--monotonic derivators.

\begin{pro}\label{pro:chain rule implicit-version}
  Let $f$ be a real\--valued function defined in a neighborhood {$U_{t^*}\subset[a,b]$ of~$t^*$ containing~$t\in [a,b]$,} $f$ $g$\--differentiable at~$t$, and  $h$ a real function defined on $f(U_{t^*})$.
  \begin{enumerate}
 \item 	 If $t^*\notin D_g$ or $f(t^{*+})=f(t^*)$, and $h$ is differentiable at $f(t^*)$, then so is $h\circ f$ and
 \[(h\circ f)_g'(t)=
 	h'\big(f(t^*)\big) f_g'(t).\]
\item If $t^*\in D_g$, $f(t^{*+})\neq f(t^*)$ and $h$ is continuous from the right at $f(t^{*+})$, then  $h\circ f$ is differentiable at $t$ and  \[(h\circ f)_g'(t)=\frac{h\big(f(t^{*+})\big)-h\big(f(t^*)\big)}{f(t^{*+})-f(t^*)}  f_g'(t).\]
\end{enumerate}
\end{pro}
\begin{proof}
We start by observing that when $h$ is differentiable at $f(t^*)$, then there exist $\gamma>0$ and a function~$\phi$ such that
$\phi(f(t^*))=0$, and
\[
h(y)=h(f(t^*))+ [h'(f(t^*))+\phi(y)][y-f(t^*)],
\]
for $y\in(f(t^*)-\gamma,f(t^*)+\gamma)\cap f(U_{t^*})$.

  We now distinguish three cases:

{\bf Case 1:} $t^*\notin D_g\cup N_g$, so $t^*=t$.  Since $f$ is $g$\--differentiable at~$t$, and $g$ is continuous at~$t$,  then, $f$ is continuous at~$t$ and  there exists $\delta_1>0$ such that
\[
|f(s)-f(t)|<\gamma,\quad\text{for }s\in (t-\delta_1, t+ \delta_1)\cap U_{t^*}.
\]
Moreover, since $f$ is $g$\--differentiable at~$t$ by Definition~\ref{dfn:g-derivative(non monotonic case)}, there exist $\delta\in(0,\delta_1]$ and a function~$\psi$ continuous at $t$  such that
$\psi(t)=0$, and
\[
f(s)=f(t)+ [f_g'(t)+\psi(s)][g(s)-g(t)],\quad\text{for all } s\in(t-\delta,t+\delta)\cap U_{t^*},\text{ with } g(s)\neq g(t).
\]
Thus, for $s\in(t-\delta,t+\delta)\cap U_{t}$,
\begin{align*}
h(f(s)) = & h(f(t))+ [h'(f(t))+\phi(f(s))][f(s)-f(t)] \\
  = & h(f(t))+ [h'(f(t))+\phi(f(s))][f_g'(t)+\psi(s)][g(s)-g(t)]\\
  = & h(f(t))+ [h'(f(t))f_g'(t)+ h'(f(t))\psi(s)+\phi(f(s))f_g'(t)+ \phi(s)\psi(s)][g(s)-g(t)]\\
  = & h(f(t))+ [h'(f(t))f_g'(t)+ \theta(s)][g(s)-g(t)],
\end{align*}
where $\theta:=h'(f(t))\cdot\psi+(\phi\circ f)\cdot f_g'(t)+ \phi\cdot\psi$. The function $\theta$ is continuous at $t$ and satisfies $\theta(t)=0$. Thus $h\circ f$ is $g$\--differentiable at $t$ and $(h\circ f)_g'(t)=h'(f(t))f_g'(t)$.

{\bf Case 2:} {$t^*\in N_g^-$}, so $t^*=t$. Since $f$ is $g$\--differentiable at~$t$, and $t\in N_g^-$,  then $f$ is continuous from the left at $t$, and there exists $\delta_1>0$ such that
\[
|f(s)-f(t)|<\gamma,\quad\text{for } s\in (t-\delta_1, t)\cap U_{t}.
\]
Since $f$ is $g$\--differentiable at~$t$, by Definition~\ref{dfn:g-derivative(non monotonic case)}, there exist $\delta\in(0,\delta_1]$ and a function $\psi$ continuous from the left at $t$  such that
$\psi(t)=0$, and
\[
f(s)=f(t)+ [f_g'(t)+\psi(s)][g(s)-g(t)],\quad\text{for all } s\in(t-\delta,t),\text{ with } \quad g(s)\neq g(t).
\]
Arguing as in the previous case, we obtain, for $s\in(t-\delta,t)$,
\begin{align*}
h(f(s)) = & h(f(t))+ [h'(f(t))+\phi(f(s))][f(s)-f(t)] \\
  = & h(f(t))+ [h'(f(t))f_g'(t)+ h'(f(t))\psi(s)+\phi(f(s))f_g'(t)+ \phi(s)\psi(s)][g(s)-g(t)]\\
  = & h(f(t))+ [h'(f(t))f_g'(t)+ \theta(s)][g(s)-g(t)],
\end{align*}
 where $\theta:=h'(f(t))\cdot\psi+(\phi\circ f)\cdot f_g'(t)+ \phi\cdot\psi$. The function $\theta$ is continuous from the left at $t$ and satisfies $\theta(t)=0$. Thus $h\circ f$ is $g$\--differentiable at $t$ and $(h\circ f)_g'(t)=h'(f(t))f_g'(t)$.

{\bf Case 3:} {$t^*\in N_g^+$. Arguing analogously to  the previous case, via replacing~$t$ by~$t^*$ and left-continuity by right-continuity,} we conclude that there exists $\delta>0$ and a function  $\theta$ continuous from the right at~{$t^*$} satisfying $\theta(t^*)=0$ such that
\begin{align*}
h(f(s)) = & h(f(t^*))+ [h'(f(t^*))+\phi(f(s))][f(s)-f(t^*)] \\
  = & h(f(t^*))+ [h'(f(t^*))f_g'(t)+ h'(f(t^*))\psi(s)+\phi(f(s))f_g'(t)+ \phi(s)\psi(s)][g(s)-g(t^*)]\\
  = & h(f(t^*))+ [h'(f(t^*))f_g'(t)+ \theta(s)][g(s)-g(t^*)],
\end{align*}
for $s\in(t^*,t^*+\delta)\cap U_{t^*}$. Thus $h\circ f$ is $g$\--differentiable at $t$ and $(h\circ f)_g'(t)=h'(f(t^*))f_g'(t)$.

{\bf Case 4:} {$t^*\in D_g$}. By Definition~\ref{dfn:g-derivative(non monotonic case)} { and Remark~\ref{rem:ld}}, since $f$ is $g$\--differentiable at~$t$, observe that
\[
f_g'(t)=\frac{f(t^{*+})-f(t^*)}{\Delta^+g(t^*)},
\]
and in particular, $f(t^{*+})$ exists. Thus, we distinguish two subcases. If $f(t^{*+})=f(t^*)$, then $f$ is continuous from the right at $t^*$,  $f_g'(t)=0$, and arguing as in {\bf Case 3} we obtain that
\[
(h\circ f)_g'(t)=h'(f(t^*))f_g'(t)=0.
\]
Now, if $f(t^{*+})\neq f(t^*)$, then $f_g'(t)\ne 0$ and, since $h$ is continuous from the right at $f(t^{*+})$, the function
\[\eta({y})=\frac{h({y})-h(f(t^*))}{{y}-f(t^*)}-\frac{h(f(t^{*+}))-h(f(t^*))}{f(t^{*+})-f(t^*)}\]
is continuous from the right at $f(t^{*+})$ and $\eta(f(t^{*+}))=0$. {For $y$ sufficiently close to $f(t^{*+})$, w}e can write
\[h({y})=h(f(t^*))+\left(\frac{h(f(t^{*+}))-h(f(t^*))}{f(t^{*+})-f(t^*)}+\eta({y})\right)({y}-f(t^*)).\]
{S}ince $f$ is $g$\--differentiable at~$t$, there exist $\delta\in \bR^+$ and a function $\psi$ continuous from the right at $t^*$  such that
$\psi(t^*)=0$, and
\[
f(s)=f(t^*)+ [f_g'(t)+\psi(s)][g(s)-g(t^*)],\quad\text{for all } s\in(t^*,t^* +\delta)\cap U_{t^*},\, \text{ with } g(s)\neq g(t^*).
\]
Thus, for $s\in(t^*,t^* +\delta)\cap U_{t^*}$ with $g(s)\neq g(t^*)$,
\begin{align*}h(f(s))=& h(f(t^*))+\left(\frac{h(f(t^{*+}))-h(f(t^*))}{f(t^{*+})-f(t^*)}+\eta(f(s))\right)(f(s)-f(t^*))\\=& h(f(t^*))+\left(\frac{h(f(t^{*+}))-h(f(t^*))}{f(t^{*+})-f(t^*)}+\eta(f(s))\right)[f_g'(t)+\psi(s)][g(s)-g(t^*)]\\= & h(f(t^*))+\left[\frac{h(f(t^{*+}))-h(f(t^*))}{f(t^{*+})-f(t^*)}f_g'(t) + \theta(s) \right] [g(s)-g(t^*)],\end{align*}
where
\[\theta(s)=\left[\frac{h(f(t^{*+}))-h(f(t^*))}{f(t^{*+})-f(t^*)}f_g'(t) + \eta(f(s)) \right]\psi(s)+\eta(f(s))f'_g(t).\]
%
The function $\theta$ is continuous from the right at $t^*$ and satisfies $\theta(t^*)=0$. Thus $h\circ f$ is $g$\--differentiable at $t$ and
\[
(h\circ f)_g'(t)=\frac{h(f(t^{*+}))-h(f(t^*))}{f(t^{*+})-f(t^*)}f_g'(t).\qedhere
\]
\end{proof}

In~\cite[Theorem 3.2]{MarSlav2024generChainRule}, Márquez Albés and Slavík provided explicit formulae to compute the Stieltjes derivative of a composition involving monotonic derivators. Building upon their argument, we present the following proposition, which provides an explicit formula for computing the Stieltjes derivative of a composition involving the non\--monotonic derivator~$g$.

\begin{pro}\label{pro:chain rule explicit-version}
  Let $f$ be a real\--valued function defined in a neighborhood {$U_{t^*}\subset[a,b]$}  of~$t^*$; containing~{$t\in[a,b]$}, and $h\in C^1(\overline{co\{f(U_{t^*})\}},\mathbb{R})$, where $\overline{co\{f(U_{t^*})\}}$ denotes the closed convex hull of $f(U_{t^*})$. If $f$ is $g$\--differentiable at~$t$, then so is $h\circ f$ and
  \[
 (h\circ f)_g'(t)= f_g'(t)\int_{0}^{1} h'\big(f(t^*)+r f_g'(t)\Delta^+g(t^*)\big)\, \operatorname{d}r.
 \]
\end{pro}
\begin{proof}
From Proposition~\ref{pro:chain rule implicit-version}, we distinguish two cases.

{\bf Case 1:} $t^*\notin D_g$. In this case  $\Delta^+ g(t^*)=0$, hence
\[
(h\circ f)_g'(t)=h'\big(f(t^*)\big) f_g'(t)=f_g'(t)\int_{0}^{1} h'\big(f(t^*)\big)\, \operatorname{d}r=f_g'(t)\int_{0}^{1} h'\big(f(t^*)+r f_g'(t)\Delta^+g(t^*)\big)\, \operatorname{d}r.
\]

{\bf Case 2:} $t^*\in D_g$. In this case, we distinguish two subcases. If $f(t^{*+})=f(t^*)$, then $f_g'(t)=0$ and hence
\[
(h\circ f)_g'(t)=h'\big(f(t^*)\big) f_g'(t)=0=f_g'(t)\int_{0}^{1} h'\big(f(t^*)+r f_g'(t)\Delta^+g(t^*)\big)\, \operatorname{d}r.
\]
Now, if $f(t^{*+})\neq f(t^*)$, then
\[
  (h\circ f)_g'(t)= \frac{h(f(t^{*+}))-h(f(t^*))}{f(t^{*+})-f(t^*)}f_g'(t) = \frac{\int_{f(t^{*})}^{f(t^{*+})}h'(\tau)\, \operatorname{d}\tau}{f(t^{*+})-f(t^*)}f_g'(t).
\]
Using the change of variable,
\[
r=\frac{\tau-f(t^*)}{f(t^{*+})- f(t^*)},
\]
and since $f_g'(t)=\frac{f(t^{*+})- f(t^*)}{\Delta^+g(t^*)}$, we obtain
\begin{align*}
  (h\circ f)_g'(t)= &f_g'(t) \int_{0}^{1}h'\big(f(t^*)+r(f(t^{*+})- f(t^*)) \big)\, \operatorname{d}r  = f_g'(t)\int_{0}^{1}h'\big(f(t^*)+r\Delta^+g(t^*)f_g'(t) \big)\, \operatorname{d}r.
\end{align*}
\end{proof}

\section{$g$\--continuity and $g$\--absolute continuity}
From now on, {let us consider $a,b\in \mathbb{R}$ with $a<b$, and  a derivator $g\in \operatorname{\operatorname{BV}^{lc}}([a,b],{\mathbb R})$. The derivator $g$ defines a pseudometric $\rho_g:[a,b]\times[a,b]\to \mathbb{R}^+$ given for $s,t \in[a,b]$ by}
\[
\rho_g(s,t)=|\widetilde{g}(s)-\widetilde{g}(t)|.
\]
The pseudometric $\rho_g$ generates a topology, which we denote $\tau_g$, over {$[a,b]$} given by its basic neighborhoods at each point {$t\in[a,b]$} by the $g$\--open balls
\[
B_g(t,r)=\{t\in \mathbb{R}: \rho_g (s,t)<r\}.
\]

In the following definition, we define the notion of $g$\--continuity with respect to $g$.
\begin{dfn}
Let {$I\subset [a,b]$}. A function $f:I \to \mathbb{R}$ is said to be {\em $g$\--continuous} at $t \in I$, if, for every $\varepsilon>0$, there exists $\delta>0$ such that
 \[
  s\in B_g(t,\delta) \Rightarrow |f(s)-f(t)|<\varepsilon.
 \]
\end{dfn}

The following proposition can be proven using the same argument as in \cite[Proposition~3.2]{FP}.
\begin{pro}
Let $f:[a,b]\to \mathbb{R}$ be a $g$\--continuous function on $[a,b]$. Then the following {statements} hold:
\begin{enumerate}
  \item $f$ is continuous from the left at each $t\in (a,b]$.
  \item If $g$ is continuous at $t\in [a,b]$, then so is $f$.
  \item If $g$ is constant on some interval $[u,v]\subset [a,b]$, then so is $f$.
\end{enumerate}
\end{pro}
\begin{rem}\label{rem:measurability of g-cont functs}
It is important to note the measurability property of a $g$\--continuous function $f:I\to \mathbb{R}$ defined on  a Borel set {$I\subset [a,b]$}. Since in this case $f:(I,\tau_g) \to(\mathbb{R},\tau_u)$ is continuous, it follows that~$f$ is Borel measurable. Using the same reasoning as in \cite[Corollary~3.5]{FP}, we can also conclude that~$f$ is Lebesgue–Stieltjes measurable.
\end{rem}

We denote by $\mathcal{BC}_{g}([a,b],{\mathbb R})$ the space of $g$\--continuous functions which are bounded on the interval $[a,b]$. Analogously to~\cite[Theorem~3.4]{FP}, the space $\mathcal{BC}_{g}([a,b],{\mathbb R})$ equipped with supremum norm
\[
\|f\|_0=\sup_{t\in[a,b]}|f(t)|,\quad\text{for all } f\in\mathcal{BC}_{g}([a,b],{\mathbb R}),
\]
is a Banach space.

\subsection{$g$\--absolute continuity}
In classical differentiation theory, it is possible to reconstruct certain functions by integrating their derivatives, provided the function is absolutely continuous. A similar concept arises in the Stieltjes derivative framework, where absolute continuity is defined with respect to a left\--continuous, non\--decreasing function $g:\mathbb{R} \to \mathbb{R}$, as described in~\cite[Definition 5.1]{PR}. In the case of derivators of {bounded} variation that are not necessarily monotonic, a corresponding notion was introduced in~\cite[Definition 6.1]{FrTo}. This definition specifically addresses derivators $g: I \to \mathbb{R}$ with \emph{controlled variation}~\cite[Definition 3.1]{FrTo} on an interval $I \subset \mathbb{R}$. The definition also incorporates the monotonic behavior of the derivator on countable connected subsets of~$I$. In~\cite{MaiaTojo-1-2025}, {this concept was extended to} derivators {$g \in\operatorname{\operatorname{BV}^{lc}}([a,b], \mathbb{R})$, $[a,b]\subset\mathbb{R}$}, without requiring them to exhibit controlled variation.

\begin{dfn}
A map $F:[a,b]\to\mathbb{R}$ is {\em $g$\--absolutely continuous}, if, for every $\varepsilon > 0$, there exists $\delta > 0$ such that, for any family
$\{(a_j , b_j)\}_{j=1}^{j=m}$ of pairwise disjoint open subintervals of $[a,b]$,
\[
\sum_{j=1}^{m} \operatorname{var}_g[a_j,b_j] <\delta \Rightarrow \sum_{i=1}^{m} |F(b_j)-F(a_j)| < \varepsilon.
\]
\end{dfn}
We denote $\mathcal{AC}_{g}([a,b],{\mathbb R})$ the set of $g$\--absolutely continuous on the interval $[a,b]$.

In particular,~\cite[Proposition~5.3]{PR} yields the following lemma.
\begin{lem}[Proposition~5.3,\cite{PR}]\label{lem:g-abs-cont func has bounded var}
 $\mathcal{AC}_{g}([a,b],{\mathbb R}) \subset\operatorname{\operatorname{BV}^{lc}}([a,b],{\mathbb R})$.
\end{lem}

Arguing analogously to~\cite[Lemma~6.8]{FrTo}, we establish conditions that guarantee the $g$\--absolute continuity of a composition of functions.
\begin{pro}\label{pro:composit of g-AbsConti functs requires lipschitz Condt}
  Let $f\in \mathcal{AC}_g([a,b],\mathbb{R})$, and $h:f([a,b]) \to\mathbb{R}$ a Lipschitz function. Then, the composition $h\circ f \in \mathcal{AC}_g([a,b],\mathbb{R})$.
\end{pro}
\begin{proof}
Let $L>0$ denote the Lipschitz constant of the function~$h$, and let us fix $\epsilon>0$. Since  $f\in \mathcal{AC}_g([a,b],\mathbb{R})$, then for $\varepsilon=\frac{\epsilon}{L}>0$, there is some $\delta>0$ such that, for any family $\{(a_j,b_j)\}_{j=1}^{m}$ of pairwise disjoint open subintervals of $[a,b]$,
\[
\sum_{j=1}^{m}\operatorname{var}_g[a_j,b_j]<\delta \quad
\implies \quad
\sum_{j=1}^m|f(b_j)-f(a_j)|<\varepsilon.
\]
Now, for a fixed family $\{(a_j,b_j)\}_{j=1}^{m}$ of pairwise disjoint open subintervals of $[a,b]$, such that $\sum_{j=1}^{m}\operatorname{var}_g[a_j,b_j]<\delta$, we have that
\[
  \sum_{j=1}^{m}|h \circ f(b_j)-h\circ f(a_j)|  \le \sum_{j=1}^{m} L|f(b_j)-f(a_j)|  < L\varepsilon = \epsilon.
\]
Hence, the composition~$h\circ f$ is $g$\--absolutely continuous on $[a,b]$.
\end{proof}
{In the following lemma, we combine the results established in~\cite[Lemma~4.7 and Theorem~5.1]{MaiaTojo-1-2025} concerning non-monotonic derivators of bounded variation}.
\begin{lem}\label{lem:LS-primitive of a func in L1g is g-abs-cont and F'g=f}
 {Let $g\in\operatorname{\operatorname{BV}^{lc}}([a,b], \mathbb{R})$ be non-constant,} and $f\in L^1_g([a,b),{\mathbb R})$, and set $F(t):=\int_{[a,t)}f \,\operatorname{d} \mu_g$. Then
 \begin{enumerate}
   \item $F\in \mathcal{AC}_{g}([a,b],{\mathbb R})$;
   \item  $F_g'=f$ $|\mu_g|$\--a.e in $[a,b]$.
 \end{enumerate}
\end{lem}

Now, we state the Fundamental Theorem of Calculus for the Lebesgue\--Stieltjes integrals involving the integral of the Stieltjes derivative, the reader is referred to~\cite[Theorem~5.2]{MaiaTojo-1-2025} for the proof.
\begin{thm}[Fundamental Theorem of Calculus for the Lebesgue\--Stieltjes integrals]\label{thm:FTC1}
 {Let $g\in\operatorname{\operatorname{BV}^{lc}}([a,b], \mathbb{R})$ be non-constant}. The following statements are equivalent:
\begin{enumerate}
 \item $F \in \mathcal{AC}_{g}([a,b],{\mathbb R})$;
 \item \begin{enumerate}
 \item $F_g'$ exists $|\mu_g|$\--a.e. in $[a,b]$;
 \item $F_g' \in L^1_g([a,b),{\mathbb R})$;
 \item $F(t)=F(a)+ \int_{[a,t)}F_g'\, \operatorname{d} \mu_g,$ for $|\mu_g|$\--a.e. $t\in [a,b]$.
 \end{enumerate}
\end{enumerate}
\end{thm}

The following result refines Lemma~\ref{lem:LS-primitive of a func in L1g is g-abs-cont and F'g=f} by assuming higher regularity, which ensures differentiability at every point. {While in the classical setting a continuous function always admits a differentiable primitive with $F'=f$ everywhere (such as in~\cite[Proposition A.2.8(iv)]{AL}), our non\--monotonic framework reveals that an extra assumption---automatically satisfied in the classical case of derivation---is required to obtain a corresponding result for the Stieltjes derivative, {as demonstrated} in~\cite[Theorem~5.3]{MaiaTojo-1-2025}}.
\begin{thm}\label{thmfcf} {Let $g\in\operatorname{\operatorname{BV}^{lc}}([a,b], \mathbb{R})$ such that $b\notin N_g^+$.}
 Let us define the function $\varphi:[a,b]\to{\mathbb R}$ by {
\begin{equation*}
\varphi(t):=\begin{dcases}
 \liminf_{s\to t}\left|\frac{g(s)-g(t)}{\widetilde{g}(s)-\widetilde{g}(t)}\right| & \mbox{if } t^*\in[a,b]\setminus(D_g \cup N_g), \\
 \liminf_{s\to t^{*+}}\left|\frac{g(s)-g(t^*)}{\widetilde{g}(s)-\widetilde{g}(t^*)}\right| & \mbox{if } t^*\in(D_g\cup N_g^+),\\
 \liminf_{s\to t^-}\left|\frac{g(s)-g(t)}{\widetilde{g}(s)-\widetilde{g}(t)}\right| & \mbox{if } t^*\in N_g^- ,
 \end{dcases}
\end{equation*}}
Assume that
\[
\varphi(t)>0\quad \text{for all } t\in[a,b].
\]
Let $f:[a,b]\to\mathbb{R}$ be a $g$\--continuous function with $f\in L^1_g([a,b),\mathbb{R})$, then the following statements hold:
\begin{enumerate}
 \item $F(t):=\int_{[a,t)} f\, \operatorname{d} \mu_g \in \mathbb{R}$ for all $t\in[a,b]$.
 \item $F \in \mathcal{AC}_{g}([a,b],{\mathbb R})$.
 \item $F_g'(t)=f(t^*)$ for all $t\in [a,b]$.
\end{enumerate}
\end{thm}

It is worthwhile to mention that the condition $\varphi(t)>0$ for all $t\in[a,b]$ occurring in Theorem~\ref{thmfcf} is clearly satisfied for each $t\in[a,b]$ such that $t^*\in D_g$. In particular,
\[
\varphi(t)=\liminf_{s\to t^{+}}\left|\frac{g(s)-g(t)}{\widetilde{g}(s)-\widetilde{g}(t)}\right|=\frac{|\Delta^+ g(t)|}{\Delta^+ \widetilde g(t)}=1>0,\text{ for all }t\in[a,b]\cap D_g.
\]
Thus, one can observe that this condition is required only on the continuity points~$t$ of~{$g$} such that {$t^*\notin D_g$}, as outlined in~\cite[Remark~5.4]{MaiaTojo-1-2025}. One can also notice that this condition naturally disappears in the {monotonic} case.

\section{Relatively compact sets of $\mathcal{BC}_{g}([a,b],\mathbb{R})$}
Let $[a,b]\subset \mathbb{R}$ be fixed. Throughout this section, we provide sufficient conditions for subsets of  $\mathcal{BC}_{g}([a,b],\mathbb{R})$ to be relatively compact with respect to topology induced by the supremum norm in the case of a general derivator {$g\in\operatorname{\operatorname{BV}^{lc}}([a,b], \mathbb{R})$.}


Since the derivator $\widetilde{g}:{[a,b]}\to\mathbb{R}$ is left\--continuous and nondecreasing, we can apply the relative compactness result {of} Frigon and López Pouso~\cite[Theorem~4.2]{FP}.  This result is based on Fra\u{n}kov\'{a}'s Theorem~\cite[Proposition~2.3]{Frank} which characterizes relative compactness in spaces of regulated functions in the case of a nondecreasing $g$. Given that $g$\--continuous functions may not be regulated at the points of~$D_g$, as demonstrated in~\cite[Example~3.3]{FP}, additional assumptions were required to control the behavior of these functions at such points.

\begin{lem}[Theorem~4.2,\cite{FP}]\label{lem:relativ-compact-set of regulated functs}
Let {$g\in\operatorname{\operatorname{BV}^{lc}}([a,b], \mathbb{R})$} and $\mathcal{D}$  a family regulated functions $f:[a,b]\to\mathbb{R}$. The family~$\mathcal{D}$ is relatively compact with respect to topology induced by the supremum norm if it satisfies the following conditions:
\begin{enumerate}
  \item $\mathcal{D}$ is $g$\--equicontinuous, i.e. for each $\epsilon>0$ and $t\in [a,b]$, there exists $\delta>0$ such that
      \[
      |f(s)-f(t)|<\epsilon \quad\text{for all $s\in [a,b]$ such that $|\widetilde{g}(s)-\widetilde{g}(t)|<\delta$, and for all $f\in\mathcal{D}$};
      \]
  \item The set $\{f(a):f\in\mathcal{D}\}$ is bounded;
  \item The family $\mathcal{D}$ has uniform right\--hand side limits in {$D_g$}, i.e. for every $t\in D_g$ and {$\epsilon>0$}, there exists $\delta>0$ such that, for each $f\in\mathcal{D}$,
      \[
      |f(s)-f(t^+)|<\epsilon\quad\text{for all $s\in(t,t+\delta)$;}
      \]
  \item The family $\mathcal{D}$ has uniformly bounded jump discontinuities in {$D_g$}, i.e. for every $t\in D_g$, there exists $\eta_t>0$ such that for each $f\in\mathcal{D}$,
      \[
      |f(t)-f(t^+)|\leq \eta_t.
      \]
\end{enumerate}
\end{lem}

{Given that} $\mathcal{AC}_g([a,b],\mathbb{R})\subset\mathcal{BC}_{g}([a,b],\mathbb{R})$, and considering the Fundamental Theorem of Calculus for Lebesgue\--Stieltjes, as established in Theorem~\ref{thm:FTC1}~\cite[Theorem~5.2]{MaiaTojo-1-2025}, we can now present sufficient conditions that ensure the relative compactness of a set in $\mathcal{AC}_g([a,b],\mathbb{R})$ within $\mathcal{BC}_{g}([a,b],\mathbb{R})$.

\begin{thm}\label{thm:relat-compact-BCg}
Let $\mathcal{D} \subset \mathcal{AC}_g([a,b],\mathbb{R})$ be a family such that $\{F(a): F\in \mathcal{D}\}$ is bounded. If there exists a function $h \in L^{1}_{g}([a,b),[0,+\infty))$ such that
$$
|F_g'(t)|\le h(t) \quad \text{ for $g$-almost all } t\in [a,b), \text{ and all } F \in \mathcal{D},
$$
then, $\mathcal{D}$ is relatively compact in $\mathcal{BC}_g([a,b],\mathbb{R})$.
\end{thm}
\begin{proof}
The argument is based on Lemma~\ref{lem:relativ-compact-set of regulated functs}. Let us prove that the family~$\mathcal{D}$ is $g$\--equicontinuous. Let $\epsilon>0$ be fixed. Since $h\in L^1_g([a,b),\mathbb{R})$, then there exists $\delta>0$ such that $\int_E h\, \operatorname{d}|\mu_g|<\epsilon$ for every $|\mu_g|$\--measurable set $E\subset [a,b)$ such that $|\mu_g|(E)<\delta$.

Let $t\in [a,b]$. Using the Fundamental Theorem of Calculus for Lebesgue\--Stieltjes integrals, Theorem~\ref{thm:FTC1}, for $s\in (t,b]$ such that $|\widetilde{g}(s)-\widetilde{g}(t)|<\delta$, we have
\[
|F(s)-F(t) |=\left|\int_{[t,s)} F_g'(\tau)\, \operatorname{d}\mu_{g}(\tau)\right|\leq\int_{[t,s)} \left|F_g'(\tau)\right|\, \operatorname{d}|\mu_{g}|(\tau)\leq \int_{[t,s)} h(\tau)\, \operatorname{d}|\mu_{g}|(\tau)<\epsilon,
\]
for all $F\in\mathcal{D}$ since $|\widetilde{g}(s)-\widetilde{g}(t)|=|\mu_g|([t,s))$. Similarly,  for $s\in [a,t)$ such that $|\widetilde{g}(s)-\widetilde{g}(t)|<\delta$, we have
\[
|F(s)-F(t) |<\epsilon, \quad\text{for all $F\in\mathcal{D}$.}
\]
In addition, the family $\mathcal{D}$ has uniform right\--hand side limits in {$D_g$}. Indeed let $t\in D_g$ and $\epsilon>0$.
Since $h\in L^1_g([a,b),\mathbb{R})$, then there exists $\delta>0$ such that $\int_{(t,s)} h(\tau)\, \operatorname{d}|\mu_{g}|(\tau)<\epsilon$ for $s\in(t,b]$ such that $|\mu_g|((t,s))<\delta$. As $\widetilde{g}(t^+)$ exists, it follows that there exists $\delta_t>0$ such that
\[
|\mu_g|((t,s))=|\widetilde{g}(s)-\widetilde{g}(t^+)|<\delta, \quad\text{for all $s\in(t,t+\delta_t)$.}
\]
Therefore, for all $s\in(t,t+\delta_t)$
\[
|F(s)-F(t^+) |=\left|\int_{(t,s)} F_g'(\tau)\, \operatorname{d}\mu_{g}(\tau)\right|\leq\int_{(t,s)} \left|F_g'(\tau)\right|\, \operatorname{d}|\mu_{g}|(\tau)\leq \int_{(t,s)} h(\tau)\, \operatorname{d}|\mu_{g}|(\tau)<\epsilon,
\]
for all $F\in\mathcal{D}$. Furthermore, the family $\mathcal{D}$ has uniformly bounded jump discontinuities in {$D_g$}, as for every $t\in D_g$,
\[
|F(t^+)-F(t)|= \left|\int_{\{t\}}F_g'(s)\, \operatorname{d}\mu_{g}(s)\right|\leq \int_{\{t\}}\left|F_g'(s)\right|\, \operatorname{d}|\mu_{g}|(s)\leq \int_{\{t\}}h(s)\, \operatorname{d}|\mu_{g}|(s),
\]
for all $F\in\mathcal{D}$. Hence, Lemma~\ref{lem:relativ-compact-set of regulated functs} ensures that the family $\mathcal{D}$ is relatively compact.
\end{proof}
\section{$g$\--exponential function}
In this section, we explore the exponential {function} within the context of linear Stieltjes differential equations involving non\--monotonic derivators. { Let $[a,b]\subset\mathbb{R}$, and $g\in\operatorname{\operatorname{BV}^{lc}}([a,b], \mathbb{R})$.}

\begin{dfn}\label{df:g-exp}
	Let $h \in L^1_g([a,b),\mathbb{R})$.
	We define the \emph{$g$-exponential function} $e_h(\cdot;a) : [a,b] \to {\mathbb{R}}$ for every $t\in [a,b]$ as
\begin{equation}\label{eq:g-exponential}
	e_h(t;a)=\prod_{s\in [a,t)\cap D_g} \left(1+h(s)\,\Delta^+ g(s)\right) \, \exp \left(\int_{[a,t)\bs D_g} h(s)\dif \mu_g \right).
\end{equation}
\end{dfn}


\begin{rem}\label{rem:til(g)-exponential func notation}
To avoid any potential ambiguity in the notation when $e_h(\cdot;a)$ is used with respect to a different derivator in the following section, we will write $e_h(\cdot;a;g)$ when it is necessary to specify which derivator is involved. Moreover, observe that if $h \in L^1_g([a,b),\mathbb{R}^+)$, and $g$ is nondecreasing, i.e. $g=\widetilde{g}$, then $e_h(t;a)=e_h(\cdot;a;\widetilde{g})\ge 1$.
\end{rem}

Before delving into the properties of the exponential function defined above, we present the following proposition which {ensures} that the exponential function is well\--defined. The proof of the second point of the proposition follows the same argument as in~\cite[Lemma~3.1]{ThesisMarquezAlbes}. We include the proof for sake of completeness.
\begin{pro}\label{pro:noreasonance conditn ensur converg log sum}
Let $h\in L^1_g([a,b),\mathbb{R})$. Then, the following statements hold:
\begin{enumerate}
  \item The set $T_h^-:=\{t\in D_g: 1+h(t) \Delta g^+(t)\le 0\}$ is finite.
  \item If $1+h(t)\Delta^+g(t) \neq 0$ for all $t\in D_g$, then
\begin{equation*}
\sum_{t \in  D_g}  \big|\log|1 + h(t)\Delta^+g(t)|\big| < \infty.
\end{equation*}
\end{enumerate}
In particular, if $1+h(t)\Delta^+g(t)>0$ for every {$t\in D_g$}, then
\begin{equation*}
\sum_{t \in  D_g}\big|\log(1 + h(t)\Delta^+g(t))\big| < \infty.
\end{equation*}
\end{pro}
\begin{proof}
1. For $t\in D_g$ such that $1+h(t)\Delta^+g(t)\le 0$, we have that $h(t)\Delta^+g(t)\le -1$, thus $|h(t)\Delta^+g(t)|\ge 1$. Therefore,
\[
\operatorname{Card} (T_h^-)=\sum_{t\in T_h^-} 1 \le \sum_{t\in T_h^-}|h(t)\Delta^+g(t)|\le \sum_{t\in D_g}|h(t)||\Delta^+g(t)| \le \int_{[a,b)}|h|\operatorname{d}|\mu_g|<\infty.
\]

2. Given that the $\log$ function satisfies $\lim\limits_{y\to 0} \left|\frac{\log|1+y|}{y}\right|=1$, this implies that there exists $\delta>0$ such that for all $|y|<\delta$, we have that
\[
\left|\Big|\frac{\log|1 +y|}{y}\Big|-1\right|<1.
\]
Thus,
\[
\left|\log|1+y|\right|<2|y|,\text{ for all } |y|<\delta.
\]
Let us consider the set $\mathfrak{D}_\delta:=\big\{t\in D_g: \delta\le |h(t)\Delta^+g(t)| \big\}$.

First, observe that $\operatorname{Card} (\mathfrak{D}_{\delta})<\infty$. Indeed, as $h \in L^1_g([a,b),\mathbb{R})$, we obtain
\[
\delta \operatorname{Card} (\mathfrak{D}_{\delta})=\sum_{t\in \mathfrak{D}_{\delta}} \delta<\sum_{t\in \mathfrak{D}_{\delta}} |h(t)\Delta^+g(t)| \le \sum_{t\in D_g} |h(t)| |\mu_g|(\{t\})\le \int_{[a,b)} |h(t)|\, \operatorname{d}|\mu_g|(t)<\infty.
\]
Thus, if we set $\mathfrak{C}_{\delta}=D_g\bs \mathfrak{D}_{\delta}$, we obtain

\begin{align*}
\sum_{t \in D_g} \left|\log |1+h(t)\Delta^+g(t)|\right| =& \sum_{t \in\mathfrak{D}_{\delta}}\left|\log |1+h(t)\Delta^+g(t)|\right| +\sum_{t \in\mathfrak{C}_{\delta}} \left|\log |1+h(t)\Delta^+g(t)|\right| \\
\le& \sum_{t \in\mathfrak{D}_{\delta}}\left|\log |1+h(t)\Delta^+g(t)|\right|+\sum_{t \in\mathfrak{C}_{\delta}}2 |h(t)\Delta^+g(t)|\\
\le& \sum_{t \in\mathfrak{D}_{\delta}}\left|\log |1+h(t)\Delta^+g(t)|\right|+\sum_{t \in D_g} 2|h(t)\Delta^+g(t)|\\
\le& \sum_{t \in\mathfrak{D}_{\delta}}\left|\log |1+h(t)\Delta^+g(t)|\right|+2\int_{[a,b)} |h(t)||\mu_g|(t)<\infty .
\end{align*}
\end{proof}

\begin{pro}
 Let $h \in L^1_g([a,b),\mathbb{R})$. The function $e_h(\cdot;a)$ defined in~\eqref{eq:g-exponential} is well\--defined.
\end{pro}

\begin{proof}
First of all, observe that since $h \in L^1_g([a,b),\mathbb{R})$, then $\big|\int_{[a,b)\bs D_g}h\operatorname{d}\mu_g\big| \le\int_{[a,b)}|h|\operatorname{d}|\mu_g|<\infty$.
Let $t\in [a,b]$. If there exists $\tau\in [a,t)\cap D_g$ such that $1+h(\tau)\Delta^+g(\tau)=0$, then
\[
e_h(t;a)=\prod_{s\in [a,t)\cap D_g} \left(1+h(s)\,\Delta^+ g(s)\right) \, \exp \left(\int_{[a,t)\bs D_g} h(s)\dif \mu_g \right)=0.
\]
Otherwise,  $1+h(\tau)\Delta^+g(\tau)\neq 0$ for all $\tau\in [a,t)\cap D_g$, so using Proposition~\ref{pro:noreasonance conditn ensur converg log sum}~(2), we deduce that
\begin{align*}
  \left|\log\big|e_h(t;a)\big|\right| = & \left|\log\left(\prod_{s\in [a,t)\cap D_g} \big|1+h(s)\,\Delta^+ g(s)\big| \, \exp \Big(\int_{[a,t)\bs D_g} h(s)\dif \mu_g(s)\Big)\right)\right| \\
  = &   \Big|\sum_{s \in[a,t)\cap D_g} \log |1+h(s)\Delta^+g(s)| + \int_{[a,t)\bs D_g} h(s)\dif \mu_g(s)\Big| \\
  \le  & \sum_{s \in[a,t)\cap D_g} \left|\log |1+h(s)\Delta^+g(s)|\right| + \int_{[a,t)\bs D_g} |h(s)|\dif |\mu_g|(s)\\
  \le  &  \sum_{s \in[a,t)\cap D_g}\big|\log |1+h(s)\Delta^+g(s)|\big| + \int_{[a,b)} |h| \dif |\mu_g|\\
    \le  &  \sum_{s \in D_g}\big|\log |1+h(s)\Delta^+g(s)|\big| + \int_{[a,b)} |h| \dif |\mu_g|\\
  <&\infty.
\end{align*}
\end{proof}

Now we introduce some notations. Let $h \in L^1_g([a,b),\mathbb{R})$ and define the sets $T_h^N:=\{ t\in D_g: 1+h(t) \Delta g^+(t)< 0\}$ and $T_h^0:=\{t\in D_g: 1+h(t)\Delta g^+(t)= 0\}$ which are subsets of the {set} $T_h^-$. Then, Proposition~\ref{pro:noreasonance conditn ensur converg log sum}~(1) ensures that $T_h^N$ and~$T_h^0$ have finite cardinality. Now, let us define
\[
\tau_0=\begin{dcases}
  \min T_h^0, & \mbox{if }  T_h^0\neq\emptyset,\\
b, & \mbox{otherwise}.
\end{dcases}
\]
Furthermore, we define
\begin{equation}\label{eq:exp:hat h}
	\overline{h}(t) = \begin{dcases}
		h(t), &\text{if $t \in [a,b)\bs D_g$,}
		\\
		\frac{\log\big|1 + h(t)\Delta^+g(t)\big|}{\Delta^+g(t)}, &\text{if $t \in [a,\tau_0) \cap D_g$,}
		\\
		0, & \text{if $t \in [\tau_0,b) \cap D_g$.}
	\end{dcases}
\end{equation}

In the next lemma, we enhance the  Stieltjes integrability of the function $\overline{h}$ defined in~\eqref{eq:exp:hat h}, this will guarantee nice properties of the $g$\--exponential function $e_h(\cdot;a)$.

\begin{lem}\label{lem:exp:hat h is g-int}
	Let $h \in L^1_g([a,b),\mathbb{R})$. Then, $\overline{h}\in L^1_g([a,b),\mathbb{R})$.
\end{lem}
\begin{proof} From Proposition~\ref{pro:noreasonance conditn ensur converg log sum}~(2), it follows that
	\begin{align*}
		\int_{[a,b)} |\overline{h}| \operatorname{d}|\mu_g|= &  \int_{D_g} |\overline{h}| \operatorname{d}|\mu_g| + \int_{[a,b)\bs D_g} \overline{h} \operatorname{d}|\mu_g| \\
		= & \sum_{t\in D_g} |\overline{h}(t)|\Delta^+g(t)  + \int_{[a,b)\bs D_g} |h|\operatorname{d}|\mu_g|  \\
		= & \sum_{t\in [a,\tau_0)\cap D_g} |\overline{h}(t)|\Delta^+g(t)  + \int_{[a,b)\bs D_g} |h|\operatorname{d}|\mu_g|  \\
		\le & \sum_{t\in  [a,\tau_0)\cap D_g} \big|\log|1+h(t)\Delta^+g(t)|\big|  + \int_{[a,b)}|h| \operatorname{d}|\mu_g|\\
		< & \infty.
	\end{align*}
\end{proof}

\begin{thm}\label{thm:g-exponential using hat h}
Let $h \in L^1_g([a,b),\mathbb{R})$.  If we set $\{t_i\}_{i=1}^{\kappa}:=\{t\in T_h^N: t<\tau_0\}$ and $t_{\kappa+1}:= \tau_0$ with $a\le t_1<\cdots<t_{\kappa}<t_{\kappa+1}=\tau_0$, then we can express the exponential function given in~\eqref{eq:g-exponential} as
\begin{equation}\label{eq:g-exponential using hat h}
e_h(t;a):=\begin{dcases}
            \exp\big(\int_{[a,t)} \overline{h}\dif\mu_g\big), & \mbox{if } a\le t\le t_1, \\
            (-1)^i \exp\big(\int_{[a,t)} \overline{h}\dif\mu_g\big), & \mbox{if } t_i<t\le t_{i+1},\,i\in\{1,\dots,\kappa\}, \\
            0, &  \mbox{if } \tau_0<t\le b.
          \end{dcases}
\end{equation}\end{thm}

\begin{proof}
  Let $t\in[a,b]$. If $t\in [a,t_1]$, then $1+h(s)\,\Delta g(s)>0$ for all $s\in [a,t)\cap D_g$, and
\begin{align*}
e_h(t;a)=&\prod_{s\in [a,t)\cap D_g} \big(1+h(s)\,\Delta^+ g(s)\big) \, \exp \left(\int_{[a,t)\bs D_g} h(s)\dif \mu_g \right) \\
= & \prod_{s\in [a,t)\cap D_g} \exp\left(\log|1+h(s)\,\Delta^+ g(s)|\right) \, \exp \left(\int_{[a,t)\bs D_g} h(s)\dif \mu_g \right)\\
  = & \exp\left(\sum_{s\in [a,t)\cap D_g} \log|1+h(s)\,\Delta^+ g(s)|\right) \, \exp \left(\int_{[a,t)\bs D_g} h(s)\dif \mu_g \right) \\
  = & \exp\left(\int_{[a,t)\cap D_g} \frac{\log|1+h(s)\,\Delta^+ g(s)|}{\Delta g(s)}\dif\mu_g(s)\right) \, \exp \left(\int_{[a,t)\bs D_g} h(s)\dif \mu_g \right) \\
  = & \exp\left(\int_{[a,t)} \overline{h}\dif\mu_g\right).
\end{align*}
To complete the proof of the other cases, we can apply the same argument piece\--by\--piece in $(t_i,t_{i+1}]$ for $i\in\{1,\dots,\kappa\}$, while it is clear that $e_h(t;a)=0$ for $t\in(\tau_0,b]$. Indeed, it suffices to observe that
\[
\prod_{s\in [a,t)\cap D_g}(1+h(s)\,\Delta^+ g(s))=(-1)^i\prod_{s\in [a,t)\cap D_g}|1+h(s)\,\Delta^+ g(s)|,
\]
for $t\in (t_i,t_{i+1}]$ and $i\in\{1,\dots,\kappa\}$.
\end{proof}

\begin{cor}
{Under} the hypotheses of Theorem~\ref{thm:g-exponential using hat h}, if $1+h(t)\Delta^+g(t) \neq 0$ for all {$t\in D_g$}, then $e_h(\cdot;a)$ never vanishes, {Moreover, if} $1+h(t)\Delta^+g(t) > 0$ for all {$t\in D_g$} then $e_h(\cdot;a)>0$.
\end{cor}

The next theorem guarantees that the exponential function can be viewed as a solution of a linear Stieltjes differential equation $|\mu_g|$\--almost everywhere, the idea of the proof is analogous to the proofs in~\cite{FP,Mar}.
\begin{thm}\label{thm:x0 exp solves linEq}
Let $h \in L^1_g([a,b),\mathbb{R})$.   Then, the function defined for every $t\in [a,b]$ by
\begin{equation*}
x(t):=x_0 e_h(t;a),
\end{equation*}
is the unique $g$\--absolutely continuous solution of  the initial value problem
\begin{equation}\label{eq:exp:linear eq}
x_g'(t) = h(t)x(t) \quad \text{for $|\mu_g|$-almost every $t \in [a,b]$}, \quad x(a) =x_0.
\end{equation}
Equivalently,
\[
x(t)= x_0+ \int_{[a,t)} h(s)x_0e_h(s;a)\, \operatorname{d}\mu_g(s) \quad \text{for every $t \in [a,b]$}.
\]
\end{thm}
\begin{proof}
Let us set $\{t_i\}_{i=1}^{\kappa}:=\{t\in T_h^N: t<\tau_0\}$ and $t_{\kappa+1}:= \tau_0$ with $a\le t_1<\cdots<t_{\kappa}<t_{\kappa+1}=\tau_0$.   From Lemmata~\ref{lem:LS-primitive of a func in L1g is g-abs-cont and F'g=f} and~\ref{lem:exp:hat h is g-int}, it follows that $\int_{[a,\cdot)} \overline{h}\, \operatorname{d}\mu_g$ is $g$\--absolutely continuous on $[a,b]$, and
  \[
  \left(\int_{[a,\cdot)} \overline{h}\, \operatorname{d}\mu_g \right)_g'(t)=\overline{h}(t),\quad \text{ for $|\mu_g|$-almost every  $t\in [a,b]$}.
  \]
From Theorem~\ref{thm:g-exponential using hat h} and Proposition~\ref{pro:composit of g-AbsConti functs requires lipschitz Condt}, it follows that $e_h(\cdot;a)$ is $g$\--absolutely continuous on $[a,b]$.  Now, it remains to show that~\eqref{eq:exp:linear eq} holds for $|\mu_g|$\--almost all ${[a,b]}\bs C_g$ combining Theorem~\ref{thm:g-exponential using hat h} and Proposition~\ref{pro:chain rule explicit-version}. Given~\eqref{eq:g-exponential using hat h}, we will prove that~\eqref{eq:exp:linear eq} holds on the set $[t_1,b]\bs C_g$, the result can be easily deduced on $[a,t_1]  \bs C_g$. Without loss of generality,  assume that $\tau_0<b$. For $|\mu_g|$\--almost every $t\in [t_1, \tau_0]\bs C_g$, we have that $t\in [t_i,t_{i+1})$ for some $i\in \{1,\dots,\kappa\}$. Thus,  by Proposition~\ref{pro:chain rule explicit-version}
\begin{equation}\label{eq:proof:deriv g-expo using chain rule}
\begin{aligned}
  \left(e_h(\cdot;a)\right)_g'(t)= & (-1)^i\,\overline{h}(t)\int_0^1 e^{\int_{[a,t)} \overline{h} \, \operatorname{d}\mu_g +r\overline{h}(t)\Delta^+g(t)}\, \operatorname{d}r\\
  =&(-1)^i\,\overline{h}(t)e^{\int_{[a,t)} \overline{h} \, \operatorname{d}\mu_g}\int_0^1 e^{r\overline{h}(t)\Delta^+g(t)}\, \operatorname{d}r= \overline{h}(t)e_h(t;a)\int_0^1 e^{r\overline{h}(t)\Delta^+g(t)}\, \operatorname{d}r.
\end{aligned}
\end{equation}
If $t\notin D_g$, then $\Delta^+g(t)=0$, $\overline{h}(t)=h(t)$, and
\begin{equation}\label{eq:proof:deriv g-expo on continuity pt}
\left(e_h(\cdot;a)\right)_g'(t)=h(t)e_h(t;a).
\end{equation}
While if $t\in D_g$, then~\eqref{eq:proof:deriv g-expo using chain rule} becomes
\begin{equation}\label{eq:proof:deriv g-expo on discontinuity pt}
\begin{aligned}
  \left(e_h(\cdot;a)\right)_g'(t)=& \frac{\log|1+h(t)\Delta^+g(t)|}{\Delta^+g(t)}e_h(t;a)\int_0^1 e^{r\log|1+h(t)\Delta^+g(t)|}\, \operatorname{d}r\\
  =& \begin{dcases}
        \frac{\log|1+h(t)\Delta^+g(t)|}{\Delta^+g(t)}e_h(t;a)\frac{e^{\log|1+h(t)\Delta^+g(t)|}-1}{\log|1+h(t)\Delta^+g(t)|}, & \mbox{if }h(t)\neq 0,\\
        0, & \mbox{otherwise,}
     \end{dcases}\\
  =& \begin{dcases}
        \frac{e_h(t;a)}{\Delta^+g(t)}h(t)\Delta^+g(t), & \mbox{if }h(t)\neq 0,\\
        0, & \mbox{otherwise,}
     \end{dcases}\\
  =& h(t)e_h(t;a).
\end{aligned}
\end{equation}
If $t=\tau_0$, then $1+h(\tau_0)\Delta^+ g(\tau_0)=0$, and we obtain
\[
 \left(e_h(\cdot;a)\right)_g'(\tau_0)=\frac{e_h(\tau_0^+;a) - e_h(\tau_0;a)}{\Delta^+ g(\tau_0)}= \frac{-e_h(\tau_0;a)}{\Delta^+ g(\tau_0)}=h(\tau_0)e_h(\tau_0;a).
\]
Now, if $t>\tau_0$, then $e_h(t;a)=0$. Thus, $\left(e_h(\cdot;a)\right)_g'(t)=0=h(t)e_h(t;a)$.

Since $e_h(\cdot;a)$ is $g$\--absolutely continuous on $[a,b]$, it follows from the Fundamental Theorem of Calculus for Lebesgue\--Stieltjes Integrals, Theorem~\ref{thm:FTC1}, and  from~\eqref{eq:proof:deriv g-expo on continuity pt} and~\eqref{eq:proof:deriv g-expo on discontinuity pt} that
\[
e_h(t,a) = e_h(a,a)+ \int_{[a,t)} h(s)e_h(s,a)\, \operatorname{d}\mu_g(s) \quad \text{for every $t \in [a,b]$}.
\]
Hence,
\[
x(t)=x_0e_h(t,a) =x_0+ \int_{[a,t)} h(s)x_0e_h(s,a)\, \operatorname{d}\mu_g(s) \quad \text{for every $t \in [a,b]$}.\qedhere
\]
\end{proof}

\section{Systems of Stieltjes differential equations with vector\--valued derivators}
We devote this section to an application of the results established in the previous sections in the theory of Stieltjes differential equations. To do so, we divided this section into three subsections. The first subsection explores continuity notions pertaining to vector\--valued derivators $\mathbf{g}=(g_1,\dots,g_n):[a,b]\to \mathbb{R}^n$, where {$g_i \in \operatorname{\operatorname{BV}^{lc}}([a,b],{\mathbb R})$} for all $i=1,\dots,n$, drawing upon the findings from~\cite{MT}.  Additionally, we present a result similar to Theorem~\ref{thmfcf} guaranteeing the Stieltjes differentiability everywhere of $g_i$\--absolutely continuous functions whose $g_i$\--derivative is not necessarily $g_i$\--continuous. In the second subsection, we consider a system of Stieltjes differential equations involving several derivators:
\begin{equation*}
\begin{aligned}
  \mathbf{x}_{\mathbf{g}}'(t)= &\mathbf{f}(t,\mathbf{x}(t)) \quad\text{for $|\mu_{\mathbf{g}}|$-almost every } t\in[t_0,t_0+T], \\
  \mathbf{x}(t_0)= & \mathbf{x}_0=(\mathbf{x}_{0,1},\dots,\mathbf{x}_{0,n}) \in\mathbb{R}^n,
\end{aligned}
\end{equation*}
where $\mathbf{f}:[t_0,t_0+T]\times \mathbb{R}^n\to \mathbb{R}^n$, {$\mathbf{g}=(g_1,\dots,g_n):[t_0,t_0+T]\to \mathbb{R}^n$, where {$g_i \in \operatorname{\operatorname{BV}^{lc}}([t_0,t_0+T],{\mathbb R})$} for all $i=1,\dots,n$} in the same fashion of~\cite{MT}. Building upon the findings from the previous sections, we establish conditions that ensure the Stieltjes differentiability of solutions throughout their entire domain. This mirrors classical results, where~$g_i$ is left\--continuous and nondecreasing, while also enhancing the intriguing conditions necessary for our more general case. In the last subsection, we provide existence and uniqueness results for local and global solutions of~\eqref{eq:SDE}.

\subsection{Continuity notions with respect to vector\--valued derivators}
Throughout this subsection, and in the spirit of~\cite{PM3,MT}, we present continuity notions relative to a vector\--valued derivator $\mathbf{g}=(g_1,\dots,g_n):[a,b]\subset \mathbb{R}\to \mathbb{R}^n$, where {$g_i \in \operatorname{\operatorname{BV}^{lc}}([a,b],{\mathbb R})$}  for all $i=1,\dots,n$. In doing so, we will be considering the vector\--valued derivator $\widetilde{\mathbf{g}}:=(\widetilde{g_1},\dots,\widetilde{g_n}):{[a,b]}\to \mathbb{R}^n$, where $\widetilde{g_i}$ is the variation function of the derivator~$g_i$ for $i=1,\dots,n$ in the sense of Definition~\ref{dfn:variation function}.

Let us denote by $\|\cdot\|$ the maximum norm in $\mathbb{R}^n$, i.e. $\|\mathbf{u}\|=\max\{|u_1|,\dots,|u_n|\}$ for $\mathbf{u}=(u_1,\dots,u_n)\in\mathbb{R}^n$.
\begin{dfn}\label{dfn:continuity notions}
Let $I\subset{[a,b]}$, and $\mathbf{f}=(f_1,\dots,f_n):I\to \mathbb{R}^n$.
\begin{enumerate}
  \item $\mathbf{f}$ is said to be {\em $\mathbf{g}$\--continuous} at $t\in I$ if, for every $i=1,\dots,n$ $f_i$ is $g_i$\--continuous.
  \item $\mathbf{f}$ is said to be {\em $\overrightarrow{\mathbf{g}}$\--continuous} at $t\in I$ if, for every $\epsilon>0$ there exists $\delta>0$ such that
\[
\text{for all $s\in I$}\quad \|\widetilde{\mathbf{g}}(s)-\widetilde{\mathbf{g}}(t)\| <\delta, \quad \|f(s) - f(t)\| < \epsilon.
\]
\end{enumerate}
\end{dfn}
We denote
\[
\mathcal{BC}_{\mathbf{g}}([a,b],\mathbb{R}^n):= \prod_{i=1}^{n}\mathcal{BC}_{g_i}([a,b],\mathbb{R}),
\]
equipped with the supremum norm $\|\cdot\|_\infty$:
\[
\|\mathbf{f}\|_\infty=\sup_{t\in[a,b]}\|\mathbf{f}(t) \| , \quad \text{ for all } \mathbf{f}\in\mathcal{BC}_{\mathbf{g}} ([a,b],\mathbb{R}^n).
\]
and
\[
\mathcal{BC}_{\overrightarrow{\mathbf{g}}}([a,b],\mathbb{R}^n) := \big\{\mathbf{f} : [a,b] \to \mathbb{R}^n : \text{$\mathbf{f}$ is $\overrightarrow{\mathbf{g}}$-continuous and bounded}\big\}.
\]
\begin{rem}\label{rem: g-continuity implies vect(g)-continuity}
  In the case where $\mathbf{g}=(g_1,\dots,g_n)$ with $g_i$ is nondecreasing and left\--continuous for all $i=1,\dots,n$, Definition~\ref{dfn:continuity notions}~(1) (resp. Definition~\ref{dfn:continuity notions}~(2)) coincides with~\cite[Defintion~3.1]{PM3} (resp.~\cite[Defintion~3.14]{MT}). In this monotonic case, the authors in~\cite[Page~18]{MT} have shown the difference between both definitions emphasizing that $\mathbf{g}$\--continuity implies~$\overrightarrow{\mathbf{g}}$\--continuity. The converse is  not always true (as shown in~\cite[Example~3.15]{MT}) except in the case when all the derivators $g_i$ generate the same topology, i.e. $C_{g_i}=C_{g_j}$ and $D_{g_i}=D_{g_j}$, for all $i,j\in\{1,\dots,n\}$,  the reader is referred to~\cite[Proposition~3.16]{MT} for more details.

   In our non\--monotonic case, since $\widetilde{\mathbf{g}}=(\widetilde{g_1},\dots,\widetilde{g_n})$ and each $\widetilde{g_i}$ is a nondecreasing and left\--continuous for $i=1,\dots,n$, it follows from the previous point that
\[
\text{$\mathbf{g}$-continuity}  \implies \text{$\overrightarrow{\mathbf{g}}$-continuity}.
\]
In particular, we have
\[
\mathcal{AC}_{\mathbf{g}}([a,b],\mathbb{R}^n)\subset\mathcal{BC}_{\mathbf{g}}([a,b],\mathbb{R}^n) \subset \mathcal{BC}_{\overrightarrow{\mathbf{g}}}([a,b],\mathbb{R}^n),
\]

where
\[
\mathcal{AC}_{\mathbf{g}}([a,b],\mathbb{R}^n):=\prod_{i=1}^n \mathcal{AC}_{g_i}([a,b],\mathbb{R}).
\]
\end{rem}
Given that the derivator $\widetilde{\mathbf{g}}=(\widetilde{g_1},\dots,\widetilde{g_n})$ and $\widetilde{g_i}$ is left\--continuous and nondecreasing of all $i=\{1,\dots,n\}$, then $\widehat{\widetilde{g}}:=\widetilde{g_1}+\cdots+\widetilde{g_n}$ is left\--continuous and nondecreasing as well, and we can state~\cite[Proposition~3.17]{MT} using $\widetilde{\mathbf{g}}$ and~$\widehat{\widetilde{g}}$.

\begin{pro}\label{pro:vect(g)continuity is hat(g)continuity}
Let $I\subset{[a,b]}$, and $\mathbf{f}=(f_1,\dots,f_n):I\to \mathbb{R}^n$. Then, $\mathbf{f}$ is $\overrightarrow{\mathbf{g}}$\--continuous at $t\in I$ if and only if $f_i$ is $\widehat{\widetilde{g}}$\--continuous at~$t$ for $i=1,\dots,n$.
\end{pro}
\begin{rem}\label{rem:measurability of vect(g)-cont funct}
The previous proposition  offers a straightforward characterization of $\widehat{\widetilde{g}}$\--continuity. Furthermore, given $\mathbf{f}=(f_1,\dots,f_n):I\subset [a,b]\to\mathbb{R}$ a $\overrightarrow{\mathbf{g}}$\--continuous function on a Borel  set $I\subset {[a,b]}$, then for all $i=1,\dots,n$, $f_i$ is $\widehat{\widetilde{g}}$\--continuous on~$I$. By means of Remark~\ref{rem:measurability of g-cont functs}, $f_i$ is Lebesgue\--Stieltjes measurable, and hence $f_i$ is $g_i$\--measurable and $\widehat{\widetilde{g}}$\--measurable.
\end{rem}

In the sequel, since we deal with the vector\--valued derivator $\widetilde{\mathbf{g}}=(\widetilde{g_1},\dots,\widetilde{g_n})$, we should exercise caution when using the notation~\eqref{eq:notation of t^*} relatively to each derivator~$g_i$, $i=1,\dots,n$. Thus, for sake of precision, we introduce the notation:
\begin{equation*}
 t^{*_i}=\begin{dcases}
 b_n, & \mbox{if $t\in(a_n,b_n)\subset C_{g_i}$, or ($t=a\notin D_{g_i}$ and $(a,b_n) \subset C_{g_i}$)},\\
  t, & \mbox{otherwise,}
 \end{dcases}
 \end{equation*}
for all $i=1,\dots,n$.

The next theorem is an adaptation of Theorem~\ref{thmfcf} in the case of several derivators. It strengthens Lemma~\ref{lem:LS-primitive of a func in L1g is g-abs-cont and F'g=f} assuming a weaker continuity condition to calculate the Stieltjes derivative on a larger domain. The proof follows an argument that is analogous to the one used in Theorem~\ref{thmfcf}. We include it here for the sake of precision.

\begin{thm}\label{thm:FTC for vect-val g}
{Assume that $b\notin \bigcup_{i=1}^n N_{g_i}^+$.} For $i=1,\dots,n$, let us define
\begin{equation}\label{eq:hat(C_i) set}
\widehat{\mathcal{C}}_i:=\{t\in[a,b]\bs D_{g_i}: t^{*_i} \in D_{\widehat{\widetilde{g}}}\},
\end{equation}
and  the function $\widehat{\varphi}_i:[a,b]\bs \widehat{\mathcal{C}}_i\to{\mathbb R}$ by {
\begin{equation}\label{eq:hat varphi_i function}
\widehat{\varphi}_i(t):=\begin{dcases}
 \liminf_{s\to t}\left|\frac{g_i(s)-g_i(t)}{\widehat{\widetilde{g}}(s)-\widehat{\widetilde{g}}(t)}\right| & \mbox{if } t^{*_i}\in[a,b]\bs(D_{g_i}\cup N_{g_i}), \\
 \liminf_{s\to (t^{*_i})^+}\left|\frac{g_i(s)-g_i(t^{*_i})}{\widehat{\widetilde{g}}(s)-\widehat{\widetilde{g}}(t^{*_i})}\right| & \mbox{if } t^{*_i}\in(D_{g_i}\cup N_{g_i}^+),\\
 \liminf_{s\to t^-}\left|\frac{g_i(s)-g_i(t)}{\widehat{\widetilde{g}}(s)-\widehat{\widetilde{g}}(t)}\right| & \mbox{if } t^{*_i}\in N_{g_i}^-.
 \end{dcases}
\end{equation}}
Assume that for $i=1,\dots,n$
\[
\widehat{\varphi}_i(t)>0\quad \text{for all $t\in[a,b]\bs \widehat{\mathcal{C}}_i$} .
\]
Let $f :[a,b]\to\mathbb{R}$ be $\widehat{\widetilde{g}}$\--continuous with $f\in L_{g_i}^1([a,b),\mathbb{R})$, $i=1,\dots,n$, then the following statements hold for all $i=1,\dots,n$:
\begin{enumerate}
 \item $f$ is $g_i$\--measurable.
 \item $F_i(t):=\int_{[a,t)} f\, \operatorname{d} \mu_{g_i} \in \mathbb{R}$ for all $t\in[a,b]$.
 \item $F_i \in \mathcal{AC}_{g_i}([a,b],{\mathbb R})$.
 \item $(F_i)_{g_i}'(t)=f(t^{*_i})$ for all $t\in [a,b]\bs \widehat{\mathcal{C}}_i$.
\end{enumerate}
\end{thm}

\begin{proof}
  1.  By Remark~\ref{rem:measurability of g-cont functs}, since $f$ is $\widehat{\widetilde{g}}$\--continuous, then $f$ is $g_i$\--measurable for $i=1,\dots,n$.

  2. Let $t\in[a,b]$, and $i\in\{1,\dots,n\}$ be fixed.
 \[
 |F_i(t)|\leqslant \int_{[a,t)}|f| \, \operatorname{d} |\mu_{g_i}|<\infty.
 \]
 Hence, $F_i(t):=\int_{[a,t)} f\, \operatorname{d} \mu_{g_i} \in \mathbb{R}$ for all $t\in[a,b]$.

3. Since  $f \in L^1_{g_i}([a,b),\mathbb{R})$, then according to Lemma~\ref{lem:LS-primitive of a func in L1g is g-abs-cont and F'g=f}, $F_i \in \mathcal{AC}_{g_i}([a,b],{\mathbb R})$.

4. Fix $i\in\{1,\dots,n\}$, $t\in [a,b]\bs \widehat{\mathcal{C}}_i$ and $\varepsilon\in{\mathbb R}^+$. We distinguish four cases:

 \textbf{Case 1:} {$t^{*_i}\in[a,b]\bs(D_{g_i} \cup N_{g_i})$. Then, $t^{*_i}=t$, and $t\in [a,b]\bs(D_{g_i}\cup C_{g_i} \cup N_{g_i})$.}
In this case, since \[
\widehat{\varphi}_i(t)=\liminf_{s\to t}\left|\frac{g_i(s)-g_i(t)}{\widehat{\widetilde{g}}(s)-\widehat{\widetilde{g}}(t)}\right|>0,
\]
 there exists $\delta_1>0$ such that, for $s\in [a,b]$ satisfying $|s-t|<\delta_1$, we have
\[
\frac{\widehat{\varphi}_i(t)}{2}< \left|\frac{g_i(s)-g_i(t)}{\widehat{\widetilde{g}}(s)-\widehat{\widetilde{g}}(t)}\right|,
\]
which implies that
\begin{equation}\label{eq:proof:case1:comparing g_i and hat(g) around a point}
|\widehat{\widetilde{g}}(s)-\widehat{\widetilde{g}}(t)|<\frac{2}{\widehat{\varphi}_i(t)} |g_i(s)-g_i(t)|.
\end{equation}
Since $f \in \mathcal{BC}_{\widehat{\widetilde{g}}}([a,b],\mathbb{R})$, then there exists $\delta_2\in{\mathbb R}^+$ such that
\[
|f(s)-f(t)|<\xi:=\frac{\varepsilon\widehat{\varphi}_i(t)}{2}, \text{ for } s \in [a,b] \text{ such that } |\widehat{\widetilde{g}}(s)-\widehat{\widetilde{g}}(t)|<\delta_2.
\]
Given that $t=t^{*_i}$, $t\notin \widehat{\mathcal{C}}_i$, then $t\not\in D_{\widehat{\widetilde{g}}}$, and $\widehat{\widetilde{g}}$ is continuous at $t$. Thus, there exists $\delta\in(0,\delta_1]$ such that, if $|s-t|<\delta$, then $|\widehat{\widetilde{g}}(s)-\widehat{\widetilde{g}}(t)|<\delta_2$. Thus, if $s\in(t-\delta,t+\delta)\cap[a,b]$, we obtain
\[
f(t)-\xi< f(s) < f(t)+\xi.
\]
Since $\mu_{g_i}=\mu_{g_i}^+ -\mu_{g_i}^-$, then for all $s\in(t-\delta,t+\delta)\cap[a,b]$ such that $s>t$ we obtain
\begin{equation}\label{eq:d1}
(f(t)-\xi)\mu_{g_i}^+([t,s))\leqslant \int_{[t,s)} f\, \operatorname{d} \mu_{g_i}^+ \leqslant (f(t)+\xi)\mu_{g_i}^+([t,s)),
\end{equation}
and
\[
(f(t)-\xi)\mu_{g_i}^-([t,s))\leqslant \int_{[t,s)} f\, \operatorname{d} \mu_{g_i}^- \leqslant (f(t)+\xi)\mu_{g_i}^-([t,s)),
\]
or equivalently,
\begin{equation}\label{eq:d2}
-(f(t)+\xi)\mu_{g_i}^-([t,s))\leqslant -\int_{[t,s)} f\, \operatorname{d} \mu_{g_i}^- \leqslant -(f(t)-\xi)\mu_{g_i}^-([t,s)),
\end{equation}
Since
\[
F_i(s)-F_i(t)=\int_{[t,s)} f\, \operatorname{d} \mu_{g_i}=\int_{[t,s)} f\, \operatorname{d} \mu_{g_i}^+-\int_{[t,s)} f\, \operatorname{d} \mu_{g_i}^-,
\]
then, adding \eqref{eq:d1} and \eqref{eq:d2},
\[ (f(t)-\xi)\mu_{g_i}^+([t,s))-(f(t)+\xi)\mu_{g_i}^-([t,s))\leqslant  F_i(s)-F_i(t) \leqslant (f(t)+\xi)\mu_{g_i}^+([t,s)) -(f(t)-\xi)\mu_{g_i}^-([t,s)),\]
that is,
\[ \mu_{g_i}([t,s))f(t)-\xi|\mu_{g_i}|([t,s))\leqslant F_i(s)-F_i(t) \leqslant \mu_{g_i}([t,s))f(t)+\xi|\mu_{g_i}|([t,s)).\]
Hence,
\[ -\xi|\mu_{g_i}|([t,s))\leqslant F_i(s)-F_i(t)-\mu_{g_i}([t,s))f(t) \leqslant\xi|\mu_{g_i}|([t,s)),\]
Or, equivalently,
\begin{equation}\label{eq:ineq for t>t0}
	|F_i(s)-F_i(t)-f(t)(g_i(s)-g_i(t))|\leqslant\xi|\mu_{g_i}|([t,s))=\xi(\widetilde{g_i}(s)-\widetilde{g_i}(t))\le \xi(\widehat{\widetilde{g}}(s)-\widehat{\widetilde{g}}(t)) .\end{equation}
Similarly, for all $s\in(t-\delta,t+\delta)\cap[a,b]$ such that $s<t$, we obtain
\begin{equation*}
\mu_{g_i}([s,t))f(t)-\xi|\mu_{g_i}|([s,t))\leqslant F_i(t)-F_i(s) \leqslant \mu_{g_i}([s,t))f(t)+\xi|\mu_{g_i}|([s,t)).
\end{equation*}
Therefore,
\[
-\mu_{g_i}([s,t))f(t)-\xi|\mu_{g_i}|([s,t))\leqslant F_i(s)-F_i(t) \leqslant -\mu_{g_i}([s,t))f(t)+\xi|\mu_{g_i}|([s,t)),\]
and, thus,
\begin{equation}\label{eq:ineq for t<t0}
	|F_i(s)-F_i(t)-f(t)(g_i(s)-g_i(t))|\leqslant\xi|\mu_{g_i}|([s,t))=\xi(\widetilde{g_i}(t)-\widetilde{g_i}(s))\le \xi(\widehat{\widetilde{g}}(t)-\widehat{\widetilde{g}}(s)).\end{equation}
Given that $\widehat{\widetilde{g}}$ is nondecreasing, using \eqref{eq:ineq for t>t0} and \eqref{eq:ineq for t<t0}, we deduce that, for all $s\in(t-\delta,t+\delta)\cap[a,b]$,
\[
|F_i(s)-F_i(t)-f(t)(g_i(s)-g_i(t))|< \xi|\widehat{\widetilde{g}}(s)-\widehat{\widetilde{g}}(t)|.
\]
Combining this inequality with~\eqref{eq:proof:case1:comparing g_i and hat(g) around a point} we obtain, for all $s\in(t-\delta,t+\delta)\cap[a,b]$,
\begin{equation}\label{eq:proof:Case1:|F_i(s)-F_i(t)-f(t)(g_i(s)-g_i(t))|<e|g_i(s)-g_i(t)|}
| F_i(s)-F_i(t)-f(t)(g_i(s)-g_i(t))| < \xi|\widehat{\widetilde{g}}(s)-\widehat{\widetilde{g}}(t)|<\varepsilon |g_i(s)-g_i(t)|.
\end{equation}
\textbf{Case 2:}  {$t^{*_i}\in  N_{g_i}^+$}. In this case, since
\[
\widehat{\varphi}_i(t)=\liminf_{s\to(t^{*_i})^+}\left|\frac{g_i(s)-g_i(t^{*_i})}{\widehat{\widetilde{g}}(s)-\widehat{\widetilde{g}}(t^{*_i})}\right|>0,
\]
there exists $\delta_1>0$ such that, for $s\in [a,b]$ satisfying $|s-t^*|<\delta_1$, we have
\[
\frac{\widehat{\varphi}_i(t)}{2}< \left|\frac{g_i(s)-g_i(t^{*_i})}{\widetilde{\widetilde{g}}(s)-\widetilde{\widetilde{g}}(t^{*_i})}\right|,
\]
which implies that
\begin{equation}\label{eq:proof:case2:comparing g and tilde(g) around a point}
	|\widehat{\widetilde{g}}(s)-\widehat{\widetilde{g}}(t^{*_i})|<\frac{2}{\widehat{\varphi}_i(t)} |g_i(s)-g_i(t^{*_i})|.
\end{equation}
Now, since $f$ is $\widehat{\widetilde{g}}$\--continuous at $t^{*_i}$, there exists $\delta_2>0$ such that
\[
|f(s)-f(t^{*_i})|<\xi:=\frac{\varepsilon\widehat{\varphi}_i(t)}{2},  \text{ for } s \in [a,b] \text{ such that } |\widehat{\widetilde g}(s)-\widehat{\widetilde g}(t^{*_i})|<\delta_2.
\]
Given that {$t^{*_i}\notin \widehat{\mathcal{C}}_i$}, $\widehat{\widetilde g}$ is continuous at $t^{*_i}$, so there exists $\delta\in(0,\delta_1]$ such that, if $|s-t^{*_i}|<\delta$, then $|\widehat{\widetilde g}(s)-\widehat{\widetilde g}(t^{*_i})|<\delta_2$. Thus, if $s\in(t^{*_i},t^{*_i}+\delta)\cap[a,b]$, we obtain,
\[
f(t^{*_i})-\xi< f(s) < f(t^{*_i})+\xi.
\]
Arguing similarly to the previous case, we deduce that, for all   $s\in(t,t+\delta)\cap[a,b]$,
\[
|F_i(s)-F_i(t^{*_i})-f(t^{*_i})(g(s)-g(t^{*_i}))|\leqslant\xi|\mu_{g_i}|([t^{*_i},s))=\varepsilon(\widehat{\widetilde{g}}(s)-\widehat{\widetilde{g}}(t^{*_i})).
\]
Using~\eqref{eq:proof:case2:comparing g and tilde(g) around a point}, we get
\begin{equation}\label{eq:proof:Case2:|F_i(s)-F_i(t)-f(t)(g_i(s)-g_i(t))|<e|g_i(s)-g_i(t)|}
| F_i(s)-F_i(t^{*_i})-f(t^{*_i})(g_i(s)-g_i(t^{*_i}))| < \xi|\widehat{\widetilde{g}}(s)-\widehat{\widetilde{g}}(t^{*_i})|<\varepsilon|g_i(s)-g_i(t^{*_i})|.
\end{equation}
\textbf{Case 3:} {$t^{*_i}\in N_{g_i}^-$, then $t^{*_i}=t$}. In this case, since
\[
\widehat{\varphi}_i(t)=\liminf_{s\to t^{-}}\left|\frac{g_i(s)-g_i(t)}{\widehat{\widetilde{g}}(s)-\widehat{\widetilde{g}}(t)}\right|>0,
\]
there exists $\delta_1>0$ such that, for $s\in [a,b]\cap (t-\delta_1,t)$ , we have
\[
\frac{\widehat{\varphi}_i(t)}{2}< \left|\frac{g_i(s)-g_i(t)}{\widehat{\widetilde{g}}(s)-\widehat{\widetilde{g}}(t)}\right|,
\]
which implies that
\begin{equation}\label{eq:proof:case3:comparing g_i and hat(g) around a point}
	|\widetilde{g}(s)-\widetilde{g}(t)|<\frac{2}{\widehat{\varphi}_i(t)} |g(s)-g(t)|.
\end{equation}
Arguing analogously to {\bf Case~2} and using~\eqref{eq:proof:case3:comparing g_i and hat(g) around a point}, we obtain that there exists $\delta\in{\mathbb R}^+$ such that, if  $s\in [a,b]\cap (t-\delta,t)$, then
\begin{equation}\label{eq:proof:Case3:|F_i(s)-F_i(t)-f(t)(g_i(s)-g_i(t))|<e|g_i(s)-g_i(t)|}
| F_i(s)-F_i(t)-f(t)(g_i(s)-g_i(t))| < \varepsilon |g_i(s)-g_i(t)|.
\end{equation}

\textbf{Case 4:} $ t^{*_i}\in D_{g_i}$. In this case, $t^{*_i} \in D_{\widehat{\widetilde{g}}}$ as well, and we have that
\[
\widehat{\varphi}_i(t)=\liminf_{s\to (t^{*_i})^+}\left|\frac{g_i(s)-g_i(t^{*_i})}{\widehat{\widetilde{g}}(s)-\widehat{\widetilde{g}}(t^{*_i})}\right|=\frac{|\Delta^+ g(t^{*_i})|}{\Delta^+\widetilde g(t^{*_i})}>0.
\]
Thus, there exists $\delta_1\in\mathbb{R}^+$ such that for $s\in (t^{*_i},t^{*_i}+\delta_1)$,
\[
\frac{\widehat{\varphi}_i(t)}{2}<\left|\frac{g_i(s)-g_i(t^{*_i})}{\widehat{\widetilde{g}}(s)-\widehat{\widetilde{g}}(t^{*_i})}\right|,
\]
which implies
\begin{equation}\label{eq:proof:case4:comparing g_i and hat(g) around a point}
\frac{\widehat{\varphi}_i(t)}{2} |\widehat{\widetilde{g}}(s)-\widehat{\widetilde{g}}(t^{*_i})|<|g_i(s)-g_i(t^{*_i})|.
\end{equation}
Now, let $M:=\max\left\{1,|f(t^{*_i})|\right\}$. Since $f$ is $g_i$\--integrable on $[a,b)$, there exists $\eta>0$ so that for $s\in(t^{*_i},b]$ such that $|\widehat{\widetilde{g}}(s)-\widehat{\widetilde{g}}((t^{*_i})^+)|=\mu_{\widehat{\widetilde{g}}}((t^{*_i},s))<\eta$, we have
\[
 \int_{(t^{*_i},s)}|f|\operatorname{d}|\mu_{g_i}|< \frac{\varepsilon\widehat{\varphi}_i(t)}{4}\Delta^+\widetilde{g_i}(t^{*_i}).
\]
Additionally, since there exists $\widehat{\widetilde g}((t^{*_i})^+)$, there is $\delta\in(0,\delta_1]$ such that for all $s\in(t^{*_i},t^{*_i}+\delta)\cap [a,b]$, we have
\[ |\widehat{\widetilde g}(s)-\widehat{\widetilde g}((t^{*_i})^+)|<\min\left\{\frac{\varepsilon \widehat{\varphi}_i(t)}{4M}\Delta^+\widetilde{g_i}(t^{*_i}),\eta\right\}.\]
So, for $s\in(t^{*_i},t^{*_i} +\delta)\cap [a,b]$,
\begin{align*}
& \left|F_i(s)-F_i(t^{*_i})-f(t^{*_i})(g_i(s)-g_i(t^{*_i}))\right|=
\left|\int_{[t^{*_i},s)}f\operatorname{d}\mu_{g_i}-f(t^{*_i})(g_i(s)-g_i(t^{*_i}))\right|\\
= &\left|\int_{\{t^{*_i}\}}f\operatorname{d}\mu_{g_i}+\int_{(t^{*_i},s)}f\operatorname{d}\mu_{g_i}-f(t^{*_i})(g_i(s)-g_i(t^{*_i}))\right|\\
= & \left|f(t^{*_i})(g_i((t^{*_i})^+)-g_i(t^{*_i}))+\int_{(t^{*_i},s)}f\operatorname{d}\mu_{g_i}-f(t^{*_i})(g(s)-g(t^{*_i}))\right|\\
=& \left|f(t^{*_i})(g_i(s)-g_i((t^{*_i})^+))+\int_{(t^{*_i},s)}f\operatorname{d}\mu_{g_i}\right|\\
\leqslant &|f(t^{*_i})||g_i(s)-g_i((t^{*_i})^+)|+ \int_{(t^{*_i},s)}|f|\operatorname{d}|\mu_{g_i}|\\
< &|f(t^{*_i})||\widetilde{g_i}(s)-\widetilde{g_i}((t^{*_i})^+)|+ \frac{\varepsilon\widehat{\varphi}_i(t)}{4}\Delta^+\widetilde{g_i}(t^{*_i})\\
\leqslant &|f(t^{*_i})|\frac{\varepsilon\widehat{\varphi}_i(t)}{4M}\Delta^+\widetilde{g_i}(t^{*_i})+\frac{\varepsilon\widehat{\varphi}_i(t)}{4}\Delta^+\widetilde{g_i}(t^{*_i})\\
\leqslant & \frac{\varepsilon\widehat{\varphi}_i(t)}{2}\Delta^+\widetilde{g_i}(t^{*_i})= \frac{\varepsilon\widehat{\varphi}_i(t)}{2}|\widetilde{g_i}((t^{*_i})^+)-\widetilde{g_i}(t^{*_i}))|\\
\leqslant & \frac{\varepsilon\widehat{\varphi}_i(t)}{2}|\widetilde{g_i}(s)-\widetilde{g_i}(t^{*_i}))|\leqslant \varepsilon\frac{\widehat{\varphi}_i(t)}{2}|\widehat{\widetilde{g}}(s)-\widehat{\widetilde{g}}(t^{*_i}))| < \varepsilon |g_i(s)-g_i(t^{*_i})|.
\end{align*}
The last inequality holds from~\eqref{eq:proof:case4:comparing g_i and hat(g) around a point}.

Now let us set
\[
h_i(s):=\begin{dcases}
 \frac{F_i(s)-F_i(t^{*_i})}{g_i(s)-g_i(t^{*_i})}-f(t^{*_i}), & \mbox{if } g_i(s)\neq g_i(t^{*_i}), \\
 0, & \mbox{otherwise}.
 \end{dcases}
\]
From~\eqref{eq:proof:Case1:|F_i(s)-F_i(t)-f(t)(g_i(s)-g_i(t))|<e|g_i(s)-g_i(t)|},~\eqref{eq:proof:Case2:|F_i(s)-F_i(t)-f(t)(g_i(s)-g_i(t))|<e|g_i(s)-g_i(t)|}, and~\eqref{eq:proof:Case3:|F_i(s)-F_i(t)-f(t)(g_i(s)-g_i(t))|<e|g_i(s)-g_i(t)|}, it results that there exists $\delta>0$ such that
\[
\begin{aligned}
&|h_i(s)|<\varepsilon, && \text{for }0<|s-t|<\delta,\ \text{with }g_i(t)\ne g_i(s),\ \text{and }t\not\in D_{g_i}\cup C_{g_i}\cup N_{g_i},\\
&|h_i(s)|<\varepsilon, && \text{for }0<s-t^{*_i}<\delta,\ \text{with }t^{*_i}\in D_{g_i}\cup N_{g_i}^+, \ \text{and } g_i(t^{*_i})\ne g_i(s),\\
&|h_i(s)|<\varepsilon, && \text{for }0<t^{*_i}-s<\delta,\ \text{with }t^{*_i} \in N_{g_i}^-,\ \text{and } g_i(t^{*_i})\ne g_i(s).
\end{aligned}
\]
Therefore, $h_i$ fulfills the assumptions of Definition~\ref{dfn:g-derivative(non monotonic case)}. Hence,
\[
(F_i)_{g_i}'(t)=f(t^{*_i})\quad \text{for all } t\in [a,b]\bs \widehat{\mathcal{C}}_i.\qedhere
\]
\end{proof}

\begin{rem}
The previous theorem allows to compute the Stieltjes derivative on the interval $[a,b]$ except at the continuity points $t$ of $g_i$ such that $t^{*_i} \in D_{\widehat{\widetilde{g}}}$. Moreover, in the particular case where $\mathbf{g}=(g,\dots,g)$, with {$g \in \operatorname{\operatorname{BV}^{lc}}([a,b],{\mathbb R})$}, observe that $\widehat{\widetilde{g}}=n\widetilde{g}$, and the $\widehat{\widetilde{g}}$\--continuity is reduced to $\widetilde{g}$\--continuity. This implies in particular that for $i,j\in \{ 1,\dots,n\}$, $D_{g_i}=D_{g_j}=D_{g}$, and $C_{g_i}=C_{g_j}=C_{g}$. Thus, according to~\cite[Corollary~3.11]{MT}, the derivators $\widetilde{g_i}=\widetilde{g}$ generate the same topology on {$[a,b]$}. Therefore, the sets $\widehat{\mathcal{C}}_i$ in~\eqref{eq:hat(C_i) set} are empty since $D_{g_i}=D_{\widehat{\widetilde{g}}}$ for all $i\in \{ 1,\dots,n\}$.  In this case, the statement of previous theorem would coincide with the statement of Theorem~\ref{thmfcf}.
\end{rem}

In the sequel, we shift our attention to other continuity notions involving  the $\mathbf{g}$\--continuity and $\overrightarrow{\mathbf{g}}$\--continuity, as the vector\--valued derivator $\widetilde{\mathbf{g}}$ confirms to the case studied in~\cite{MT}.
\begin{dfn}
Let $I\times X \subset[a,b]\times\mathbb{R}^n$, $\mathbf{f}: I\times X\to \mathbb{R}^n$, and $(t,\mathbf{x})\in I\times X$ be fixed.
\begin{enumerate}
  \item $\mathbf{f}$ is said to be {\em $(\mathbf{g}\times \operatorname{id}_{\mathbb{R}^n})$\--continuous at} $(t,\mathbf{x})$ if, {for $i=1,\dots,n$, and} for every $\epsilon > 0$, there exists $\delta > 0$ such that
\[
 |f_i (s,\mathbf{y}) - f_i(t,\mathbf{x})| < \epsilon,  \quad \text{for all  $(s,\mathbf{y}) \in I\times X$ such that $|\widetilde{g_i}(s)-\widetilde{g_i}(t)| <\delta$ and $\|\mathbf{x}-\mathbf{y}\|<\delta$};
\]
  \item $\mathbf{f}$ is said to be {\em $(\overrightarrow{\mathbf{g}}\times \operatorname{id}_{\mathbb{R}^n})$\--continuous at} $(t,\mathbf{x})$ if, for every $\epsilon>0$ there exists $\delta>0$ such that
\[
\|\mathbf{f}(s,\mathbf{y}) - \mathbf{f}(t,\mathbf{x})\| < \epsilon, \text{ for all $s\in I$ such that $\|\widetilde{\mathbf{g}}(s)-\widetilde{\mathbf{g}}(t)\| <\delta$ and $\|\mathbf{x}-\mathbf{y}\|<\delta$}.
\]
\end{enumerate}
\end{dfn}

\begin{pro}[{\cite[Proposition 3.22]{MT}}]
Let $I\times X \subset[a,b]\times\mathbb{R}^n$, $\mathbf{f}: I\times X\to \mathbb{R}^n$, and $(t,\mathbf{x})\in I\times X$ be fixed. If $\mathbf{f}$ is $(\mathbf{g}\times \operatorname{id}_{\mathbb{R}^n})$\--continuous at $(t,\mathbf{x})$, then $\mathbf{f}$ is $(\overrightarrow{\mathbf{g}}\times \operatorname{id}_{\mathbb{R}^n})$\--continuous at $(t,\mathbf{x})$.
\end{pro}
\begin{lem}\label{lem:(vec(g).id)-conti implies f(.x(.))vec(g)-conti}
  Let $I\times X \subset[a,b]\times\mathbb{R}^n$, $\mathbf{f}: I\times X\to \mathbb{R}^n$ be $(\overrightarrow{\mathbf{g}}\times \operatorname{id}_{\mathbb{R}^n})$\--continuous on $I\times X $. If $\mathbf{x}:I\to X$ is $\overrightarrow{\mathbf{g}}$\--continuous on $I$, then the composition $\mathbf{f}(\cdot,\mathbf{x}(\cdot))$ is $\overrightarrow{\mathbf{g}}$\--continuous on $I$.
\end{lem}
\begin{rem}\label{rem:(vec(g).id)-conti implies f(.x(.))vec(g)-conti but not til(g)-cont}
In light of Remark~\ref{rem: g-continuity implies vect(g)-continuity}, the previous lemma ensures that the composition of a $(\mathbf{g}\times \operatorname{id}_{\mathbb{R}^n})$\--continuous function~$\mathbf{f}$ with a $\mathbf{g}$\--continuous function~$\mathbf{x}$ is $\overrightarrow{\mathbf{g}}$\--continuous,  however, this is generally not sufficient to ensure the $\mathbf{g}$\--continuity of the composition $\mathbf{f}(\cdot,\mathbf{x}(\cdot))$, the reader is referred to~\cite[Remark~3.24]{MT} for a counter\--example. Nevertheless, using Proposition~\ref{pro:vect(g)continuity is hat(g)continuity}, we can use the conditions of Lemma~\ref{lem:(vec(g).id)-conti implies f(.x(.))vec(g)-conti} to exploit $\mathbf{f}(\cdot,\mathbf{x}(\cdot))$ as a $\overrightarrow{\mathbf{g}}$\--continuous function, i.e. $f_i(\cdot,\mathbf{x}(\cdot))$ as a $\widehat{\widetilde{g}}$\--continuous function for all $i=1,\dots,n$.
\end{rem}
\subsection{Everywhere solutions}
In this subsection, we consider a system  of  Stieltjes differential equations of the form
\begin{equation}\label{eq:SDE}
\begin{aligned}
  \mathbf{x}_{\mathbf{g}}'(t)= &\mathbf{f}(t,\mathbf{x}(t)) \quad\text{for $|\mu_{\mathbf{g}}|$-almost every } t\in[t_0,t_0+T], \\
  \mathbf{x}(t_0)= & \mathbf{x}_0=(\mathbf{x}_{0,1},\dots,\mathbf{x}_{0,n}) \in\mathbb{R}^n,
\end{aligned}
\end{equation}
to be understood for all $i=1,\dots,n$ as
\begin{equation*}
\begin{aligned}
  (x_i)_{g_i}'(t)= &f_i(t,\mathbf{x}(t)) \quad\text{for $|\mu_{g_i}|$-almost every } t\in[t_0,t_0+T], \\
  x_i(t_0)= & \mathbf{x}_{0,i}\in\mathbb{R},
\end{aligned}
\end{equation*}
where { $\mathbf{g}=(g_1,\dots,g_n):[t_0,t_0+T]\to \mathbb{R}$ is a vector\--valued derivator such that $g_i \in \operatorname{\operatorname{BV}^{lc}}([t_0,t_0+T],{\mathbb R})$ for all $i=1,\dots,n$, $\mathbf{f}=(f_1,\dots,f_n):[t_0,t_0+T]\times \mathbb{R}^n\to \mathbb{R}^n$. In the sequel, we provide conditions ensuring that the solutions of the problem~\eqref{eq:SDE} are Stieltjes differentiable everywhere on their definition domain.}

We start by the definition of solutions of problem~\eqref{eq:SDE}.
\begin{dfn}
A function  $\mathbf{x}=(x_1,\dots,x_n):[t_0,t_0+T] \to \mathbb{R}^n$ is said to be a {\em solution of the problem~\eqref{eq:SDE}}, if $\mathbf{x}\in\mathcal{AC}_{\mathbf{g}}([t_0,t_0+T],\mathbb{R}^n)$, $\mathbf{x}(t_0)=\mathbf{x}_0$ and
\[
(x_i)_{g_i}'(t)=f_i(t,\mathbf{x}(t))  \quad \text{ for $|\mu_{g_i}|$-almost all } t\in [t_0,t_0+T].
\]
\end{dfn}

Now, we show that the result in this section imply the corresponding results in the classical sense of ``everywhere'' $g$\--differentiable solutions making use of the $\mathbf{g}$\--continuity and $\widehat{\widetilde{\mathbf{g}}}$\--continuity conditions on the interval of existence $[t_0,t_0+T]$. In the sequel, we assume that $t_0+T \notin \bigcup_{i=1}^n  N_{g_i}^+$. In doing so, we require additional assumptions on the non\--monotonic derivator~$\mathbf{g}$:
\begin{equation*}
  ({\rm H}_{\mathbf{g}}): \qquad\qquad \prod_{i=1}^{n} \varphi_i(t)>0, \text{ for all } t\in [t_0,t_0+T],i=1\dots,n,
\end{equation*}

\begin{equation*}
  ({\rm H}_{\widehat{\mathbf{g}}}): \qquad\qquad \prod_{i=1}^{n} \widehat{\varphi}_i(t)>0, \text{ for all } t\in [t_0,t_0+T]\bs \widehat{\mathcal{C}}_i,i=1\dots,n,
\end{equation*}
where $\varphi_i:[t_0,t_0+T]\to \mathbb{R}^+$ is defined for $i=1,\dots,n$ by {
\begin{equation*}
\varphi_i(t):=\begin{dcases}
 \liminf_{s\to t}\left|\frac{g_i(s)-g_i(t)}{\widetilde{g}_i(s)-\widetilde{g}_i(t)}\right| & \mbox{if } t^{*_i}\in[t_0,t_0+T]\bs(D_{g_i}\cup N_{g_i}), \\
 \liminf_{s\to (t^{*_i})^+}\left|\frac{g_i(s)-g_i(t^{*_i})}{\widetilde{g}_i(s)-\widetilde{g}_i(t^{*_i})}\right| & \mbox{if } t^{*_i}\in(D_{g_i}\cup N_{g_i}^+),\\
 \liminf_{s\to t^-}\left|\frac{g_i(s)-g_i(t)}{\widetilde{g}_i(s)-\widetilde{g}_i(t)}\right| & \mbox{if } t\in N_{g_i}^-,
 \end{dcases}
\end{equation*}}
and $\widehat{\varphi}_i:[t_0,t_0+T]\bs \widehat{\mathcal{C}}_i\to \mathbb{R}^+$ as defined in~\eqref{eq:hat varphi_i function}.

\begin{pro}\label{pro:everywhere solution with f(.,x(.)) til g continuous}
Let $\mathbf{f}:[t_0, t_0 +T]\times\mathbb{R}^n\to\mathbb{R}^n$, and $\mathbf{x}=(x_1,\dots,x_n) : [t_0, t_0 +\tau] \to \mathbb{R}^n$ be a solution of problem~\eqref{eq:SDE}. If there exists $j\in\{1,\dots,n\}$ such that $\phi_j(t)>0$ for all $t\in [t_0, t_0 +T]$ and $f_j(\cdot,\mathbf{x}(\cdot))$ is $g_j$\--continuous on $[t_0, t_0 +T]$, then
\[
(x_j)_{g_j}'(t) = f_j(t^{*_j},\mathbf{x}(t^{*_j})) \quad \text{ for all } t\in [t_0, t_0 +T],
\]
and $(x_j)_{g_j}'$ is $g_j$\--continuous on $[t_0, t_0 +T]$. In particular, if~$(\rm{H}_{\mathbf{g}})$ holds and $\mathbf{f}(\cdot,\mathbf{x}(\cdot))$ is {$\mathbf{g}$}\--continuous on $[t_0,t_0+T]$, then
\[
(x_i)_{g_i}'(t) = f_i(t^{*_i}, \mathbf{x}(t^{*_i})) \quad \text{ for all  $t\in [t_0, t_0 +T]$ and $i=1,\dots,n$},
\]
and $(x_i)_{g_i}'$ is $g_i$\--continuous on {$[t_0, t_0 +T]$} for all $i=1,\dots,n$.
\end{pro}
\begin{proof}
Since $\mathbf{x}$ is a solution of problem~\eqref{eq:SDE}, then $\mathbf{x}\in\mathcal{AC}_{\mathbf{g}}([t_0,t_0+\tau],\mathbb{R}^n)$. Using the Fundamental Theorem of Calculus for Lebesgue\--Stieltjes integrals, Theorem~\ref{thm:FTC1}, we have that
\[
x_j(t)=\mathbf{x}_{0,j}+\int_{[t_0,t_0+t)} f_j(s,\mathbf{x}(s)) \operatorname{d}\mu_{g_j}(s),\quad \text{for all  $t\in [t_0, t_0 +T]$},
\]
and $f_j(\cdot,\mathbf{x}(\cdot))\in L^1_{g_j}([t_0, t_0 +T),\mathbb{R})$. Since  $f_j(\cdot,\mathbf{x}(\cdot))$ is $g_j$\--continuous on $[t_0, t_0 + T]$ and {$\phi_j(t)>0$}, it follows from Theorem~\ref{thmfcf} that
\[
(x_j)_{g_j}'(t)=f_j(t^{*_j},\mathbf{x}(t^{*_j})) \quad\text{for all }t\in [t_0,t_0+T].
\]
Hence, $(x_j)_{g_j}'$ is $g_j$\--continuous on {$[t_0, t_0 +T]$}.

In particular, if~$(\rm{H}_{\mathbf{g}})$ holds and $\mathbf{f}(\cdot,\mathbf{x}(\cdot))$ is $\mathbf{g}$\--continuous on $[t_0, t_0 + T]$, then $f_i(\cdot,\mathbf{x}(\cdot))$ is $g_i$\--continuous on $[t_0, t_0 + T]$, for $i=1,\dots,n$, and the result follows from the same argument for every $i\in\{1,\dots,n\}$.
\end{proof}

In light of Lemma~\ref{lem:(vec(g).id)-conti implies f(.x(.))vec(g)-conti}, and Remark~\ref{rem:(vec(g).id)-conti implies f(.x(.))vec(g)-conti but not til(g)-cont}, we can deduce that when we compose a $(\mathbf{g}\times \operatorname{id}_{\mathbb{R}^n})$\--continuous function~$\mathbf{f}$ with a solution~$\mathbf{x}$ of the problem~\eqref{eq:SDE}, the composition $\mathbf{f}(\cdot,\mathbf{x}(\cdot))$ is not necessarily~$\mathbf{g}$\--continuous. Though, we can explore some properties of ``everywhere" solutions involving a weaker continuity assumption; the $\overrightarrow{\mathbf{g}}$\--continuity. The next proposition is a generalization of~\cite[Proposition~4.4]{MT}.

\begin{pro}\label{pro: f(.,x(.)) vect(g)-cont ensures solutions 'everywhere'}
Let
 $\mathbf{f}:[t_0, t_0 +T]\times\mathbb{R}^n\to\mathbb{R}^n$, and $\mathbf{x}=(x_1,\dots,x_n) : [t_0, t_0 +T] \to \mathbb{R}^n$ be a solution of the problem~\eqref{eq:SDE}. Assume that~$(\rm{H}_{\widehat{\mathbf{g}}})$ holds. If $\mathbf{f}(\cdot,\mathbf{x}(\cdot))$ is $\overrightarrow{\mathbf{g}}$\--continuous on
$[t_0, t_0 +T]$, then
\[
(x_i)_{g_i}'(t) = f_i(t^{*_i}, \mathbf{x}(t^{*_i})) \quad \text{ for all  $t\in [t_0, t_0 +T]\bs\widehat{\mathcal{C}}_i$ and $i=1,\dots,n$},
\]
where $\widehat{\mathcal{C}}_i$ as in~\eqref{eq:hat(C_i) set}. In particular, $(x_i)_{g_i}'$ is $\widehat{\widetilde{g}}$\--continuous on $[t_0, t_0 +T]\bs\widehat{\mathcal{C}}_i$.
\end{pro}
\begin{proof}
Let $i\in\{1,\dots,n\}$ be fixed. Since $\mathbf{f}(\cdot,\mathbf{x}(\cdot))$ is $\overrightarrow{\mathbf{g}}$\--continuous, then $f_i(\cdot,\mathbf{x}(\cdot))$ is $\widehat{\widetilde{g}}$\--continuous. Since $\mathbf{x}=(x_1,\dots,x_n)$ is a solution of the problem~\eqref{eq:SDE} on $[t_0,t_0+T]$, then we have
 \[
x_i(t)=\mathbf{x}_{0,i}+\int_{[t_0,t_0+t)} f_i(s, \mathbf{x}(s))\, \operatorname{d} \mu_{g_i}(s), \quad \text{for all }t\in[t_0,t_0+T],
 \]
and $f_i(\cdot,\mathbf{x}(\cdot))\in L^1_{g_i}([t_0,t_0+T),\mathbb{R})$. Since  $f_i(\cdot,\mathbf{x}(\cdot))$ is {$\widehat{\widetilde{g}}$\--continuous} on $[t_0, t_0 + T]$ and~$(\rm{H}_{\widehat{{\mathbf{g}}}})$ holds, the result follows immediately from Theorem~\ref{thm:FTC for vect-val g}, and we obtain
\[
(x_i)_{g_i}'(t) = f_i(t^{*_i}, \mathbf{x}(t^{*_i})),\quad\text{ for all $t\in [t_0,t_0+T]\bs\widehat{\mathcal{C}}_i$}.
\]
 Hence, $(x_i)_{g_i}'$ is $\widehat{\widetilde{g}}$\--continuous on $[t_0, t_0 +T]\bs\widehat{\mathcal{C}}_i$.
\end{proof}

\begin{pro}\label{pro:everywhere solution with f (til g. id_R^n) continuous}
Let $\mathbf{f}:[t_0, t_0 +T]\times\mathbb{R}^n\to\mathbb{R}^n$ be $(\overrightarrow{\mathbf{g}}\times \operatorname{id}_{\mathbb{R}^n})$\--continuous on $[t_0, t_0 +T]$, and $\mathbf{x}=(x_1,\dots,x_n) : [t_0, t_0 +T] \to \mathbb{R}^n$ be a solution of the problem~\eqref{eq:SDE} on the interval $[t_0,t_0+T]$. Assume that~$(\rm{H}_{\widehat{\mathbf{g}}})$ holds.
Then,
\[
(x_i)_{g_i}'(t) = f_i(t^{*_i}, \mathbf{x}(t^{*_i})) \quad \text{ for all  $t\in [t_0, t_0 +T]\bs\widehat{\mathcal{C}}_i$ and $i=1,\dots,n$},
\]
where $\widehat{\mathcal{C}}_i$ as in~\eqref{eq:hat(C_i) set}. In particular, $(x_i)_{g_i}'$ is $\widehat{\widetilde{g}}$\--continuous on $[t_0, t_0 +T]\bs\widehat{\mathcal{C}}_i$.
\end{pro}
\begin{proof}  Since $\mathbf{f}$ is $(\overrightarrow{\mathbf{g}}\times \operatorname{id}_{\mathbb{R}^n})$\--continuous, and $\mathbf{x}\in \mathcal{AC}_{\mathbf{g}}([t_0, t_0 +T],\mathbb{R}^n)$, which is in particular $\overrightarrow{\mathbf{g}}$\--continuous, then Lemma~\ref{lem:(vec(g).id)-conti implies f(.x(.))vec(g)-conti} ensures that~$\mathbf{f}(\cdot,\mathbf{x}(\cdot))$ is $\overrightarrow{\mathbf{g}}$\--continuous on $[t_0, t_0 + T]$. Therefore, the result follows immediately from {Proposition~\ref{pro: f(.,x(.)) vect(g)-cont ensures solutions 'everywhere'}}.
\end{proof}
\subsection{Picard and Peano type existence results}
In this subsection, we provide existence and uniqueness results for the system of Stieltjes differential equations~\eqref{eq:SDE}. {Let $\mathbf{g}=(g_1,\dots,g_n):[t_0,t_0+T]\subset \mathbb{R}\to \mathbb{R}^n$, where $g_i \in \operatorname{\operatorname{BV}^{lc}}([t_0,t_0+T],{\mathbb R})$  for all $i=1,\dots,n$. Without loss of generality, we assume that $g_i$ is continuous at $t_0$ for all $i\in\{1,\dots,n\}$}, the reader is referred to~\cite[Sec.~5, p.~21]{FP}. We start by the definition of $\mathbf{g}$\--Carathéodory functions.
\begin{dfn}[$\mathbf{g}$\--Carath\'{e}odory function]\label{df:g-Caratheodory-R}
  Let $X\subset\mathbb{R}^n$. A function $\mathbf{f}=(f_1,\dots,f_n):[t_0,t_0+T]\times X\to \mathbb{R}^n$ is said to be {\it$\mathbf{g}$\--Carath\'{e}odory} if, for every $i\in\{1,\dots,n\}$ $f_i$ is $g_i$\--Carathéodory, i.e.  the following conditions hold:
  \begin{enumerate}
    \item for each $\mathbf{u}\in X$, $t\mapsto f_i(t,\mathbf{u})$ is $g_i$\--measurable;
    \item for $g_i$\--almost every $t\in [t_0,t_0+T]$, $\mathbf{u}\mapsto f_i(t,\mathbf{u})$ is continuous on $X$;
    \item for every $r>0$, there is a function $h_{i,r} \in L^{1}_{g_i}([t_0,t_0+T),[0,+\infty))$ such that
        \[
        |f_i(t,\mathbf{u})| \le h_{i,r}(t), \quad\text{ for $g_i$-almost all } t\in [t_0,t_0+T),  \text{ for all } \mathbf{u}\in X\cap \overline{B(\mathbf{0}_{\mathbb{R}^n},r)}.
        \]
  \end{enumerate}
\end{dfn}

Arguing as in~\cite[Lemma~7.2]{FP}, it can be shown that the composition $\mathbf{f}(\cdot, \mathbf{x}(\cdot))$ is Lebesgue\--Stieltjes integrable for every bounded $\mathbf{g}$\--continuous function~$\mathbf{x}$.
\begin{lem}\label{lem:fox-g-int}
  Let $X\subset\mathbb{R}^n$. Let $\mathbf{f}:[t_0,t_0+T]\times X\to \mathbb{R}^n$  be a $\mathbf{g}$\--Carath\'{e}odory function. Then, the map $\mathbf{f}(\cdot,\mathbf{x}(\cdot))\in L^{1}_{\mathbf{g}}([t_0,t_0+T),\mathbb{R}^n):=\prod_{i=1}^{n}L^1_{g_i}([t_0,t_0+T),\mathbb{R})$ for all $\mathbf{x}\in \mathcal{BC}_{\mathbf{g}}([t_0,t_0+T],X)$.
\end{lem}

We first establish a global existence and uniqueness result.
\begin{thm}
  Let $\mathbf{f}=(f_1,\dots,f_n):[t_0,t_0+T]\times\mathbb{R}^n\to\mathbb{R}^n$ be a function satisfying:
   \begin{enumerate}
    \item for all $i\in\{1,\dots,n\}$, $ f_i(\cdot,\mathbf{u})$ is $g_i$\--measurable for all $\mathbf{u}\in \mathbb{R}^n$;
    \item for all $i\in\{1,\dots,n\}$, $ f_i(\cdot,\mathbf{u}_0)\in L^{1}_{g_i}([t_0,t_0+T),\mathbb{R})$ for some $\mathbf{u}_0\in \mathbb{R}^n$;
    \item there is a function $L \in L^{1}_{g_i}([t_0,t_0+T),[0,+\infty))$ for all $i\in\{1,\dots,n\}$ such that
        \[
        |f_i(t,\mathbf{u})-f_i(t,\mathbf{v})| \le L(t)\|\mathbf{u}-\mathbf{v}\|,
        \]for $g_i$\--almost all $t\in [t_0,t_0+T)$ and all $\mathbf{u},\mathbf{v} \in  \mathbb{R}^n$.
  \end{enumerate}
  Then, problem~\eqref{eq:SDE} has a unique solution $\mathbf{x}\in\mathcal{AC}_{g}([t_0,t_0+T],\mathbb{R}^n)$.
\end{thm}

\begin{proof}
Let us define the operator $\mathcal{F}=(F_1,\dots,F_n):\mathcal{BC}_{\mathbf{g}}([t_0,t_0+T],\mathbb{R}^n)\to \mathcal{BC}_{\mathbf{g}}([t_0,t_0+T],\mathbb{R}^n)$ for all $\mathbf{x}\in \mathcal{BC}_{\mathbf{g}}([t_0,t_0+T],\mathbb{R}^n)$ and  $t\in [t_0,t_0+T]$ by
\[
\mathcal{F}\mathbf{x}(t)=\mathbf{x}_0+\int_{[t_0,t)}\mathbf{f}(s,\mathbf{x}(s))\, \operatorname{d}\mu_{\mathbf{g}}(s),
\]
to be understood for $i\in\{1,\dots,n\}$ as
\[
F_i\mathbf{x}(t)=\mathbf{x}_{0,i}+\int_{[t_0,t)}f_i(s,\mathbf{x}(s))\, \operatorname{d}\mu_{g_i}(s).
\]
Observe that the assumptions imply that~$\mathbf{f}$ is~$\mathbf{g}$-Carathéodory, thus, the operator $\mathcal{F}$ is well\--defined.
Let $L$ be the function given by Assumption~(3). It follows from~\cite[Proposition~2.4]{MT}, that $L\in L_{\widehat{\widetilde{g}}}([t_0,t_0+T),\mathbb{R})$. Hence, the $\widehat{\widetilde{g}}$\--exponential map $e_L(\cdot;t_0;\widehat{\widetilde{g}})$ defined with respect to the derivator $\widehat{\widetilde{g}}:=\widetilde{g_1}+\cdots+\widetilde{g_n}$ is well\--defined, as denoted in Remark~\ref{rem:til(g)-exponential func notation}. Now, let us defined the norm
\[
\|\mathbf{x}\|_L=\sup_{t\in [t_0,t_0+T]} \frac{\|\mathbf{x}(t)\|}{e_L(t;t_0;\widehat{\widetilde{g}})},\quad\text{for all }\mathbf{x}\in \mathcal{BC}_{\mathbf{g}}([t_0,t_0+T],\mathbb{R}^n).
\]
Observe that the norm $\|\cdot\|_L$ is equivalent to the norm $\|\cdot\|_{\infty}$. Indeed, since $L\geq0$, and $\widehat{\widetilde{g}}$ is nondecreasing, then by means of Remark~\ref{rem:til(g)-exponential func notation}, $e_L(t;t_0;\widehat{\widetilde{g}})\geq 1$ for all $t\in [t_0,t_0+T]$. Thus, for all $\mathbf{x}\in \mathcal{BC}_{\mathbf{g}}([t_0,t_0+T],\mathbb{R}^n)$.
\[
\|\mathbf{x}\|_L\leq \|\mathbf{x}\|_{\infty}\leq e_L(t_0+T;t_0;\widehat{\widetilde{g}}) \|\mathbf{x}\|_L .
\]
The space $\left(\mathcal{BC}_{\mathbf{g}}([t_0,t_0+T],\mathbb{R}^n),\|\cdot\|_L\right)$ is a Banach space. Thus, it is sufficient to prove that the operator $\mathcal{F}$ is a contraction so as to obtain a unique fixed point. Let $\mathbf{u},\mathbf{v} \in \mathcal{BC}_{\mathbf{g}}([t_0,t_0+T],\mathbb{R}^n)$:
\begin{align*}
  \|\mathcal{F}\mathbf{u}-\mathcal{F}\mathbf{v}\|_{L} \leq & \sup_{t\in [t_0,t_0+T]}e_L^{-1}(t;t_0;\widehat{\widetilde{g}}) \left\|\int_{[t_0,t_0+t)} \mathbf{f}(s,\mathbf{u}(s))-\mathbf{f}(s,\mathbf{v}(s))\, \operatorname{d}\mu_{\mathbf{g}}(s)\right\| \\
    = &   \sup_{t\in [t_0,t_0+T]}e_L^{-1}(t;t_0;\widehat{\widetilde{g}})\left|\int_{[t_0,t_0+t)} f_j(s,\mathbf{u}(s))-f_j(s,\mathbf{v}(s))\, \operatorname{d}\mu_{g_j}(s)\right|,
\end{align*}
for some $j\in\{1,\dots,n\}$. Thus,
\begin{align*}
  \|\mathcal{F}\mathbf{u}-\mathcal{F}\mathbf{v}\|_{L} \leq & \sup_{t\in [t_0,t_0+T]}e_L^{-1}(t;t_0;\widehat{\widetilde{g}})\int_{[t_0,t)} \Big|f_j(s,\mathbf{u}(s))-f_j(s,\mathbf{v}(s))\Big|\, \operatorname{d}|\mu_{g_j}|(s)\\
  \leq & \sup_{t\in [t_0,t_0+T]}e_L^{-1}(t;t_0;\widehat{\widetilde{g}})\int_{[t_0,t)} L(s)\|\mathbf{u}(s)-\mathbf{v}(s)\|\, \operatorname{d}|\mu_{g_j}|(s)\\
  \leq & \sum_{i=1}^n  \sup_{t\in [t_0,t_0+T]}e_L^{-1}(t;t_0;\widehat{\widetilde{g}})\int_{[t_0,t)} L(s)\|\mathbf{u}(s)-\mathbf{v}(s)\|\, \operatorname{d}|\mu_{g_i}|(s)\\
   =   & \sup_{t\in [t_0,t_0+T]}e_L^{-1}(t;t_0;\widehat{\widetilde{g}})\int_{[t_0,t)} L(s)\|\mathbf{u}(s)-\mathbf{v}(s)\|\, \operatorname{d}\mu_{\widehat{\widetilde{g}}}(s)\\
  \leq & \|\mathbf{u}-\mathbf{v}\|_L \sup_{t\in [t_0,t_0+T]}e_L^{-1}(t;t_0;\widehat{\widetilde{g}})\int_{[t_0,t)} L(s)e_L(s;t_0;\widehat{\widetilde{g}})\, \operatorname{d}\mu_{\widehat{\widetilde{g}}}(s)\\
  \leq & \|\mathbf{u}-\mathbf{v}\|_L \sup_{t\in [t_0,t_0+T]}e_L^{-1}(t;t_0;\widehat{\widetilde{g}})\Big(e_L(t;t_0;\widehat{\widetilde{g}})-1\Big)\\
  =& \Big(1-e_L^{-1}(t_0+T;t_0;\widehat{\widetilde{g}})\Big)\|\mathbf{u}-\mathbf{v}\|_L.
\end{align*}
Since $1-e_L^{-1}(t_0+T;t_0;\widehat{\widetilde{g}})<1$, then $\mathcal{F}$ is a contraction. By the Banach Contraction Principle~\cite{GD}, it follows that $\mathcal{F}$ has a unique fixed point. Therefore, problem~\eqref{eq:SDE} has a unique solution $\mathbf{x}:[t_0,t_0+T]\to \mathbb{R}^n$.
\end{proof}

In the following theorem, we provide a local existence result. Let $\overline{B(\mathbf{x}_0,r)}$ denote the open ball of $\mathbb{R}^n$  centered at $\mathbf{x}_0\in\mathbb{R}^n$ with radius $r>0$.
\begin{thm}\label{thm:unique local existence}
Let $\mathbf{f}=(f_1,\dots,f_n):[t_0,t_0+T]\times\overline{B(\mathbf{x}_0,r)}\to\mathbb{R}^n$ be a function satisfying:
   \begin{enumerate}
    \item for every $i\in \{1,\dots,n\}$, $ f_i(\cdot,\mathbf{u})$ is $g_i$\--measurable for all $\mathbf{u}\in \overline{B(\mathbf{x}_0,r)}$;
    \item for every $i\in \{1,\dots,n\}$, $ f_i(\cdot,\mathbf{u}_0)\in L^{1}_{g_i}([t_0,t_0+T),\mathbb{R})$ for some $\mathbf{u}_0\in \overline{B(\mathbf{x}_0,r)}$;
    \item for every $r>0$, there is a function $L_{r} \in L^{1}_{g_i}([t_0,t_0+T),[0,+\infty))$ for all $i\in \{1,\dots,n\}$ such that
        \[
        |f_i(t,\mathbf{u})-f_i(t,\mathbf{v})| \le L_{r}(t)\|\mathbf{u}-\mathbf{v}\|,
        \]for $g_i$\--almost all $t\in [t_0,t_0+T)$ and all $\mathbf{u},\mathbf{v} \in  \overline{B(\mathbf{x}_0,r)}$.
  \end{enumerate}
  Then, there exists $\tau\in(0,T]$ such that problem~\eqref{eq:SDE} has a unique solution $\mathbf{x}\in\mathcal{AC}_{\mathbf{g}}([t_0,t_0+\tau],\mathbb{R}^n)$.
\end{thm}
\begin{proof}
For all $i=1,\dots,n$, observe that there exists $H_i\in L^1_{g_i}([t_0,t_0+T),[0,\infty))$  such that
\begin{equation}\label{eq:proof:f(t,x(t))<H(t)}
|f_i(t,\mathbf{x})|\leq |f_i(t,\mathbf{u}_0)|+L(t)\|\mathbf{x}-\mathbf{u}_0\|\leq H_i(t):= |f_i(t,\mathbf{u}_0)|+L(t)2r,
\end{equation}
for $|\mu_g|$\--almost all $t\in [t_0,t_0+T)$ and all $\mathbf{x}\in\overline{B(\mathbf{x}_0,r)}$.

Now, let us fix $\tau\in(0,T]$ such that
\begin{equation}\label{eq:proof:int H <r}
\max_{i\in\{1,\dots,n\}}\int_{[t_0,t_0+\tau)} H_i(s)\, \operatorname{d}|\mu_{g_i}|(s)\leq r.\end{equation}
Let us define the subset
\[
\mathcal{U}:=\{\mathbf{x}\in \mathcal{BC}_{\mathbf{g}}([t_0,t_0+\tau],\mathbb{R}^n): \|\mathbf{x}(t)-\mathbf{x}_0\|\leq r \text{ for all } t\in [t_0,t_0+\tau] \},
\]
and $\mathcal{F}:\mathcal{U}\to \mathcal{U}$ by
\[
\mathcal{F}(\mathbf{x})(t)=\mathbf{x}_0+\int_{[t_0,t)} \mathbf{f}(s,\mathbf{x}(s))\, \operatorname{d}\mu_{\mathbf{g}}(s)\quad\text{ for all $t\in [t_0,t_0+\tau]$}.
\]
The operator $\mathcal{F}$ is well\--defined and maps the closed subset $\mathcal{U}$ to~$\mathcal{U}$. Indeed, by~\eqref{eq:proof:f(t,x(t))<H(t)}, and~\eqref{eq:proof:int H <r} it follows that $\mathbf{f}(\cdot,\mathbf{x}(\cdot))\in L^1_{\mathbf{g}}([t_0,t_0+\tau),\mathbb{R}^n)$ for all $\mathbf{x}\in \mathcal{U}$. Moreover, for all $\mathbf{x}\in \mathcal{U}$ and $t\in [t_0,t_0+\tau]$:
\begin{equation*}
 \|\mathcal{F}\mathbf{x}(t) - \mathbf{x}_0\|=  \left\|\int_{[t_0,t)}\mathbf{f}(s,\mathbf{x}(s))\, \operatorname{d}\mu_{\mathbf{g}}(s)\right\| = \left|\int_{[t_0,t)} f_j(s,\mathbf{x}(s))\, \operatorname{d}\mu_{g_j}(s)\right|,
\end{equation*}
for some $j\in\{1,\dots,n\}$. Therefore, by~\eqref{eq:proof:f(t,x(t))<H(t)}, we obtain
\begin{align*}
 \|\mathcal{F}\mathbf{x}(t) - \mathbf{x}_0\|\leq \int_{[t_0,t)}\left| f_j(s,\mathbf{x}(s))\right|\, \operatorname{d}|\mu_{g_j}|(s)  \leq &
 \int_{[t_0,t)} H_j(s)\, \operatorname{d}|\mu_{g_j}|(s)\leq
 \int_{[t_0,t_0+\tau)} H_j(s)\, \operatorname{d}|\mu_{g_j}|(s)\\
 \leq & \max_{i\in\{1,\dots,n\}}\int_{[t_0,t_0+\tau)} H_i(t)\, \operatorname{d}|\mu_{g_i}|(t)\leq r.
\end{align*}

Now that $\mathcal{F}:\mathcal{U}\to\mathcal{U}$ is well\--defined, arguing as in the proof of the previous theorem, if we endow the space $\mathcal{BC}_{\mathbf{g}}([t_0,t_0+\tau],\mathbb{R}^n)$ with the norm
\[
\|\mathbf{x}\|_L=\sup_{t\in[t_0,t_0+\tau]}e_L^{-1}(t;t_0;\widehat{\widetilde{g}}) \|\mathbf{x}(t)\|,
\]
 as in the proof of the previous theorem, we can show that $\mathcal{F}$ is a contraction. The conclusion follows from the Banach contraction principle.
\end{proof}
In~\cite[Theorem~8.3]{FrTo}, the authors have established a Peano\--type existence result for the problem~\eqref{eq:SDE} with nonmonotonic derivators with {\em controlled variation}. The proof uses the Schauder fixed point theorem. To demonstrate the compactness of the fixed point operator therein, the authors employed the notions of {\em $g$\--equicontinuity} and {\em $g$\--stability}~\cite[Definitions~4.10 and~4.11]{FrTo} relative to a left\--continuous derivator~$g$ with {\em controlled variation}~\cite[Definition~3.1]{FrTo}. These definitions hold true in our more general case, and we obtain the same existence result for our setting. We include a proof based on Theorem~\ref{thm:relat-compact-BCg}
for the sake of completeness.
\begin{thm}\label{thm:local existence}
Let $r >0$, and  $ \mathbf{f} :[t_0, t_0 + T] \times \overline{B(\mathbf{x}_0, r)}\to \mathbb{R}^n$ be a $\mathbf{g}$\--Carathéodory function. Then there exists $\tau \in (0,T]$ such that  problem~\eqref{eq:SDE} has at least one solution $\mathbf{x}\in\mathcal{AC}_{\mathbf{g}}([t_0,t_0+\tau],\mathbb{R}^n)$.
\end{thm}
\begin{proof}
Let $R:=\|\mathbf{x}_0\|+r$. Since the function~$\mathbf{f}$ is $\mathbf{g}$\--Carathéodory, then for each $i\in \{1,\dots,n\}$, there exists a function $h_{R,i}\in L^1_{g_i}([t_0,t_0+T),\mathbb{R})$  such that
 \[
 |f_i(t,\mathbf{x})|\leq h_{R,i}(t), \quad \text{for $|\mu_{g_i}|$-almost all $t\in [t_0,t_0+T)$ and all $\mathbf{x}\in \overline{B(\mathbf{x}_0, r)}$.}
 \]
 Now, let us fix $\tau \in (0,T]$ such that
 \[
\max_{i\in\{1,\dots,n\}}\int_{[t_0,t_0+\tau)} h_{R,i}(t)\, \operatorname{d}|\mu_{g_i}|(t)\leq r.
 \]
 Let us define $\mathcal{U}$, and $\mathcal{F}:\mathcal{U}\to \mathcal{U}$ as in the proof of the previous theorem. $\mathcal{U}$ is nonempty closed convex subset of $\mathcal{BC}_{\mathbf{g}}([t_0,t_0+\tau],\mathbb{R}^n)$. It remains to show that $\mathcal{F}=(F_1,\dots,F_n)$ is a continuous and compact operator.

 $\diamond$ Continuity: Let $\{\mathbf{x}_n\}_{n\in\mathbb{N}}\subset \mathcal{U}$ be a sequence converging to $\mathbf{x} \in \mathcal{U}$, then $\{\mathbf{x}_n\}_{n\in\mathbb{N}}$ is bounded by~$R$. For $i\in\{1,\dots,n\}$, since $f_i(t,\cdot)$ is a continuous function for $|\mu_{g_i}|$\--almost every $t\in[t_0,t_0+T]$, it follows that
\[
f_i(t,\mathbf{x}_n(t)) \to f_i(t,\mathbf{x}(t)) \quad \text{for $|\mu_{g_i}|$-almost every } t\in [t_0,t_0+T).
\]
Using Definition~\ref{df:g-Caratheodory-R}~(iii), there is a function $h_{R,i} \in L^1_{g_i}([t_0,t_0+T),[0,+\infty))$ such that
\[
|f_i(t,\mathbf{x}_n(t))|\le h_{R,i}(t) \quad \text{ for $|\mu_{g_i}|$-almost every } t\in [t_0,t_0+T).
\]
By means of Lebesgue's Dominated Convergence Theorem, we obtain
\[
\int_{[t_0,t)}f_i(s,\mathbf{x}_n(s))\, \operatorname{d}|\mu_{g_i}|(s) \ \to\  \int_{[t_0,t)}f_i(s,\mathbf{x}(s))\, \operatorname{d}|\mu_{g_i}|(s)\quad \text{ for every } t\in [t_0,t_0+T].
\]
Therefore, for some $j\in\{1,\dots,n\}$, we have that
\begin{align*}
\|\mathcal{F}(\mathbf{x}_n)-\mathcal{F}(\mathbf{x})\|_{\infty}= & \sup_{t\in[t_0,t_0+\tau]}\|\mathcal{F}\mathbf{x}_n(t)-\mathcal{F}\mathbf{x}(t)\|\\ = &
\sup_{t\in[t_0,t_0+\tau]}\left|\int_{[t_0,t)}f_j(s,\mathbf{x}_n(s))-f_j(s,\mathbf{x}(s))\, \operatorname{d}\mu_{g_j}(s)\right|\\
\leq &\sup_{t\in[t_0,t_0+\tau]}\int_{[t_0,t)}\left|f_j(s,\mathbf{x}_n(s))-f_j(s,\mathbf{x}(s))\right|\, \operatorname{d}|\mu_{g_j}|(s).
\end{align*}
Therefore,
\[
 \|\mathcal{F}(\mathbf{x}_n)-\mathcal{F}(\mathbf{x})\|_{\infty}\to 0.
\]
Hence, the operator $\mathcal{F}$ is continuous.

$\diamond$ Compactness: For every $\mathbf{x}\in \mathcal{U}$, since $\mathbf{f}$ is $\mathbf{g}$\--Carathéodory, it follows from Lemma~\ref{lem:fox-g-int} that $\mathbf{f}(\cdot,\mathbf{x}(\cdot)) \in L^{1}_{\mathbf{g}}([a,b),\mathbb{R}^n)$. Now, Lemma~\ref{lem:LS-primitive of a func in L1g is g-abs-cont and F'g=f} yields that $\mathcal{F}\mathbf{x} \in \mathcal{AC}_{\mathbf{g}}([t_0,t_0+\tau],\mathbb{R}^n)$.

Moreover, for $i\in\{1,\dots,n\}$, there is a function $h_{R,i} \in L^{1}_{g_i}([t_0,t_0+\tau),[0,+\infty))$ such that
\[
\big|(F_i\mathbf{x})_{g_i}'(t)\big|=|f_i(t,\mathbf{x}(t))|\leq h_{R,i}(t),
\]
for $|\mu_{g_i}|$\--almost all $t\in[t_0,t_0+\tau)$ and all $\mathbf{x}\in \mathcal{U}$. Since $\{\mathcal{F}\mathbf{x}(t_0):\mathbf{x}\in \mathcal{U}\}$ is bounded, it follows from Theorem~\ref{thm:relat-compact-BCg} that $\{\mathcal{F}\mathbf{x}:\mathbf{x}\in \mathcal{U}\}$ is relatively compact in $\mathcal{BC}_{\mathbf{g}}([t_0,t_0+\tau),\mathbb{R}^n)$, i.e. $\mathcal{F}(\mathcal{U})$ is relatively compact. Hence the operator $\mathcal{F}$ is compact. By the Schauder fixed point theorem~\cite{GD}, it follows that $\mathcal{F}$ has at least a fixed point in $\mathcal{U}$. Hence, problem~\eqref{eq:SDE} has at least one solution $\mathbf{x}:[t_0,t_0+\tau]\to\mathbb{R}^n$.
 \end{proof}

\section{Application to a Stieltjes-based model of solar panel and their battery performance under thermal stress}

Solar energy systems are inherently subject to variability due to the intermittent nature of sunlight. The efficiency of solar panels fluctuates throughout the day,  and energy storage systems, such as batteries, are essential for balancing supply and demand. Over time, battery health gradually deteriorates, resulting in a reduction of storage capacity. Classical approaches to modelling battery ageing and coupled PV–storage systems typically use ordinary differential equations (ODEs), representing calendar ageing with Arrhenius\--type laws, state-of-charge (SoC) weighted instantaneous decay, or cycle\--counting damage accumulators; cumulative ageing thereby arises implicitly via time integration of these local rates~\cite{BARRE2013,Plett,YaoLeiTang}.

In this work, we propose a mathematical model describing the accumulated thermal effects on photovoltaic (PV) panels and battery health, depending on solar irradiance and temperature. The system is modeled using a system of Stieltjes differential equations to capture time-dependent variations in energy storage, the degradation of battery health, and the cumulative thermal stress on the PV panel. The Stieltjes formulation differs from these ODE-based approaches by separating the "clock" that accumulates exposure from the instantaneous response: each state evolves with respect to its own derivator~\cite{PM3,MT}, which explicitly encodes cumulative effects such as degree-hours or thermal exposure.

We incorporate the effects of battery degradation over time on the system's storage capacity. {The model serves as a demonstration of the utility of the Stieltjes differential equations framework in capturing heterogeneous  time scales in coupled energy systems. The ability to explicitly model damage accumulation and recovery via non-uniform and non-monotonic clocks provides a mathematically rigorous framework to explore alternatives to classical ODE approaches. In this context, we consider a $T$-day simulation window for solar irradiance and temperature, which can support scenario analysis.}

The dynamics will be modeled using the following state variables:

\begin{itemize}
  \item $E(t)$: Energy stored in the battery at time~$t$, where $E(t) \in [0, E_{\text{max}}]$, and $E_{\text{max}}$ is the nominal battery capacity.
  \item $H(t)\in (0, 1]$: State of Battery health at time~$t$, with $H(t) = 1$ representing perfect health, and $H(t) \to 0$ representing complete degradation.
  \item $S(t)\in [0, 1]$: Cumulative thermal stress experienced by the PV panels at time~$t$, which accumulates as the system operates.  $S(t) = 0$ represents the perfect state of the panel, and $S(t)=1$ corresponds to the panel approaching its maximum cumulative thermal stress, i.e., near the limit of thermal fatigue it can tolerate.
\end{itemize}


Each state variable will be associated to a derivator such that:
\begin{itemize}
  \item $g_1$
   is the derivator of the stored energy function $E(t)$;
  \item { $g_2$} is the derivator of the battery health function $H(t)$;
  \item {$g_3$} is the derivator of the thermal stress index $S(t)$.
\end{itemize}

\paragraph*{Energy storage function $E$:}
The energy stored in the battery $E(t)$ evolves based on the dynamic balance between the power output of the solar panel and the demand function $D$ which is {an integrable} function. The model introduces a state variable $S(t)$ that tracks the cumulative thermal stress on the solar cells within the PV panel. This stress, which accumulates over time, impacts the efficiency of the entire panel's instantaneous power output~$P(t,S(t))$ from the PV panel can be given by:
\[
P(t,S(t)) := A\alpha(t,S(t)) \cdot I_{\rm POA}(t),
\]
where:
\begin{itemize}
  \item $A$: Area of the PV panel.
  \item $I_{POA}(t)$: plane-of-array (POA) irradiance incident on the module surface at time~$t$.
  \item $\alpha$: denotes the panel efficiency function. Over time, accumulated stress reduces panel efficiency. Classical models account for the effect of the cell temperature~$T_{\rm cell}(\cdot)$~\cite{KUMARLAHA2022,SunAsaDeetKia}. Therefore, an accurate thermal model is essential to improve efficiency and predict the module's lifespan. In our work, we consider $\alpha$ as a function of the PV cell temperature~$T_{\rm cell}(t)$ with a feedback from the thermal stress $S(t)$ at each time~$t$:
\[
\alpha(t,S(t)) := \alpha_{\text{ref}} \left(1 - \gamma \left( T_{\rm cell}(t) - T_{\text{ref}} \right) \right) (1-\rho S(t)).
\]
  \item $\alpha_{\text{ref}}$: reference efficiency at a reference temperature~$T_{\text{ref}}$.
  \item $\gamma$: efficiency temperature coefficient.
  \item $\rho > 0$: constant that quantifies performance degradation due to stress.
  \item $T_{\rm cell}(t)$: temperature of the cell at time~$t$, which depends on $T_{\rm ambient}(t)$, the ambient temperature, and the POA irradiance. Several models in the literature also account for wind speed and other environmental factors, see~\cite{Chapter23DWB2013,FaimanDavid2008,HasNashMor2024} and references therein. In the present work, we adopt the classical Ross relation~\cite{Ross1980}:
\[
T_{\rm cell}(t) = T_{\text{ambient}}(t) + \frac{{\rm NOCT}-20^\circ C}{800 \,W/m^2} \cdot I_{POA}(t).
\]
\item $\rm NOCT$: the module's nominal operating cell temperature (NOCT)~\cite{Chapter23DWB2013}, i.e. the cell or module temperature that is reached when the cells are mounted in their normal way at a solar radiation level of 800 $W/m^2$, a wind speed of 1$ m/s$, an ambient temperature of $20\,^\circ C$, and no-load operation;
\item $\frac{{\rm NOCT}-20^\circ C}{800\, W/m^2}$: is the thermal coefficient of the PV panel, representing its ability to convert irradiance into heat.

{This  simplified, steady-state thermal model is adopted for methodological demonstration in this application section. Using this model ensures that $T_{\rm cell}(t)$ is a known function of time, which is a necessary condition for the subsequent definition of the derivator $g_3$ to be state-independent. For applications involving rapidly changing weather conditions, where the transient effects of the PV panel's thermal inertia become significant, a dynamic thermal model—governed by its own differential equation to model $T_{\rm cell}$ as an additional time-varying component—would be required, for which the reader is referred to~\cite[Table~1]{SunAsaDeetKia}. We use the steady-state relation here for numerical simulations under a summer clear-sky scenario, where transient effects are minimal.}
\end{itemize}

Storage is an internal process that responds to both surplus and deficit conditions in real time. The efficiency of this storage process depends on several factors, including conversion losses and the current charge level of the battery. As the battery nears full capacity, the rate at which energy can be stored decreases. Moreover, the maximum storage capacity itself is modulated by the battery health factor $H(t)$ at each time~$t$. As the battery health factor decays over time, the maximum storage capacity decreases as well. Since the battery health $H(t)$ scales the effective storage capacity via $E_{\rm max}(t):=E_{\rm max} H(t)$, as the battery degrades, not only does its capacity decrease, but its self\--discharge rate increases as well. This phenomenon, whereby energy is lost over time due to internal chemical reactions, can be modeled by a health\--dependent leakage term of the form $\lambda(t)=\lambda_0(1+\delta(1-H(t))^2)$, leading to an additional decay term $-\lambda(t)E(t)$ in the storage equation, where $\delta>0$ is a modulation factor that adjusts the leakage rate $\lambda_0$ based on the state of health of the battery.

The energy storage dynamics can be described by the following equation:
\begin{equation}\label{eq:E'(t)}
\begin{dcases}
E'(t)=\frac{\eta_0H(t)(P(t,S(t))-D(t))}{1+\delta_T(T_{\rm ambient}(t)-T_{\rm op})^2}\left(1-\frac{E(t)}{E_{\rm max}H(t)}\right)-\lambda_0(1+\delta(1-H(t))^2)E(t),&\\
E(0)=E_0\in[0,E_{\rm max}H(0)]\subset[0,E_{\rm max}],&
\end{dcases}
\end{equation}
for a.e. $t\in [0,24T]$, where:
\begin{itemize}
  \item The term $\frac{\eta_0 H(t)}{1+\delta_T(T_{\rm ambient}(t)-T_{\rm op})^2}$ defines a storage efficiency function that accounts for both battery health and ambient temperature conditions, where $\eta_0\in(0,1)$ is the nominal storage efficiency under the optimal temperature $T_{\rm op}$ for battery performance, and $\delta_T>0$ is a sensitivity parameter reflecting how efficiency declines with temperature change.
  \item The term $P(t,S(t))-D(t)$ governs whether the battery is charging or discharging:
\begin{itemize}
  \item Surplus energy: If $P(t,S(t))>D(t)$, the system has more power than it needs for consumption. This surplus is available for storage in the battery.
  \item Energy deficit: If $P(t,S(t))<D(t)$, the system needs more power than it is generating. In this case, the battery will discharge to supply the difference between generated and consumed power.
\end{itemize}
  \item the term $1-\frac{E(t)}{E_{\rm max}H(t)}$ represents a scaling factor that decreases as the battery storage approaches its maximum capacity. When the battery is close to full, this factor becomes smaller, meaning less energy can be added.
\end{itemize}

\paragraph*{Battery health function $H$:}
The battery health function $H(t) \in (0,1]$ models the degradation of the battery's effective storage capacity over time, where the maximal usable capacity is given by $E_{\max} H(t)$ at each time~$t$. Battery degradation results from a combination of usage patterns and environmental stressors, notably temperature and state of charge. Empirical evidence~\cite{BARCELLONA2020,SALDANA2022} shows that the battery ages faster when frequently operated near extreme states of charge (close to 0\% or 100\%) and when exposed to temperatures exceeding a threshold $T_{\rm{thresh}} \in [30^\circ C, 40^\circ C]$. These factors often interact: thermal stress exacerbates degradation especially when the battery is at high or low charge levels, highlighting the importance of managing both variables concurrently.

We model these effects by incorporating an ambient temperature\--dependent derivator~$g_2$ and nonlinear sensitivity to the state of charge. The derivator $g_2$ of the battery health would capture the effect of prolonged thermal exposure above~$T_{\rm thresh}$ on the battery. The idea is that below~$T_{\rm thresh}$, thermal stress is negligible, and above it, aging accelerates exponentially. To this aim, the slopes of $g_2$  represent the thermal degradation factor over time. When the ambient temperature is below~$T_{\rm thresh}$, thermal stress is negligible, thus, the slopes of~$g_2$ should indicate the baseline state of the battery's degradation factor with respect to temperature. This means that no additional degradation occurs beyond the basic degradation from the charge/discharge cycles, and $g_2$ can be considered as the classical derivator during such periods. Now, if cumulative thermal exposure accumulated during some periods when the temperatures exceed the threshold~$T_{\rm thresh}$, then $g_2$ must present greater slopes during these periods. Thus, a possible choice is given for each time~$t\in [0,24T]$ by
\[
g_2(t)=\int_{0}^{t}\exp\Big(\beta_{\rm thermal}\max(0,T_{\rm ambient}(s)-T_{\rm thresh})\Big)\operatorname{d}s,
\]
where $\beta_{\rm thermal}$ is a scaling factor that controls the sensitivity of the battery's health to high temperatures above the threshold~$T_{\rm thresh}$.

Taking all these factors into consideration, the degradation of the battery health can be described by the Stieltjes differential equation
\begin{equation}\label{eq:H'g(t)}
\begin{dcases}
H_{g_2}'(t)=-\nu\left(\frac{E(t)}{E_{\rm max}H(t)}-50\%\right)^4H(t),\quad \text{for $g_2$-almost every $t\in[0,24T]$},&\\
H(0)=H_0\in(0,1],&
\end{dcases}
\end{equation}
where
\begin{itemize}
  \item $\nu>0$  determines the degradation rate of the battery health due to the charge and discharge cycles when the ambient temperature is bellow the threshold temperature~$T_{\rm thresh}$.
  \item the factor $\left(\frac{E(t)}{E_{\rm max}H(t)}-50\%\right)^4$ models the empirical observation that battery health degrades faster when the system frequently operates near full or empty states of charge. The power 4 is chosen to reflect a nonlinear sensitivity to the state of charge $\frac{E(t)}{E_{\rm max}H(t)}$ (with respect to the battery's currently degraded capacity due to aging), with minimal degradation near the optimal 50\% level, and sharply increasing degradation near the boundaries (0\% and 100\%).
\end{itemize}

\paragraph*{Thermal stress of the PV panel:}
To account for the impact of operating conditions on the PV panel, we introduce the cumulative thermal stress state variable $S(t)$. This variable encodes the accumulated stress as the panel operates under non-ideal temperatures. While extreme events such as hotspots or electrical faults can, in principle, cause sudden jumps in thermal stress, in typical operational scenarios the continuous accumulation of stress is the dominant factor. Therefore, for planning, efficiency analysis, and lifetime estimation, we focus on the continuous part. {The dynamics are can be given in terms of a derivator $g_3:[0,24T]\to\mathbb{R}$, which quantifies the instantaneous thermal stress rate based on deviations of the cell temperature $T_{\rm cell}(t)$ from the optimal operating temperature~$T_{\rm op}$ at each time~$t\in [0,24T]$. As a possible choice, we define $g_3$ for every $t\in [0,24T]$ by
\[
g_3(t) = \int_0^t \left[ \mu_1 \cdot \max(0, T_{\rm cell}(s) - T_{\text{op}})^\beta - \mu_2 \cdot \max(0, T_{\rm op} - T_{\rm cell}(s))^{\beta_r} \right] \dif s,
\]}
where:
\begin{itemize}
    \item $\mu_1$ (resp. $\mu_2$) is the thermal stress accumulation rate due to overheating  (resp. overcooling). $\mu_2$~often satisfies $\mu_2<\mu_1$which indicates that stress builds faster than it recovers;
    \item $\beta \ge 1$ controls the nonlinearity of stress accumulation above $T_{\rm op}$;
    \item $\beta_r \ge 1$ controls the nonlinearity of recovery below $T_{\rm op}$.
\end{itemize}
Observe that when the panel's temperature maintains the optimal temperature~$T_{\rm op}$, then $g_3$ becomes stationary, which implies that no stress will be accumulated.

Thus, the cumulative thermal stress on the PV panel can be modeled by:
\begin{equation}\label{eq:S'g(t)}
\begin{dcases}
    S_{g_3}'(t) = \begin{dcases}
(1-S(t)),  & \mbox{if } T_{\rm cell}(t)\geq T_{\text{op}} \\
                   S(t), & \mbox{otherwise},
                  \end{dcases}\\
    S(0)=S_0\in[0,1],
\end{dcases}
\end{equation}
for $g_3$-almost all $t\in[0,24T]$. The factors $1-S(t)$ during stress accumulation ($T \ge T_{\rm op}$) and $S(t)$ during stress recovery ($T < T_{\rm op}$) encode a state-dependent rate: the closer the panel is to its stress bounds, the slower the change. Specifically, a nearly unstressed panel accumulates stress faster and recovers slowly, while a highly stressed panel accumulates stress more slowly but relaxes faster. 

Now that we use the Stieltjes differential equations~\eqref{eq:E'(t)},~\eqref{eq:H'g(t)}, and~\eqref{eq:S'g(t)} to describe the dynamics of our system, let $R>0$ be such that $\mathcal{H}\nsubseteq\overline{B(\mathbf{x}_0,R)}$ where $\mathcal{H}:=\{(E,H,S)\in\overline{B(\mathbf{x}_0,R)}: H=0\}$, and define the functions $f_i:[0,24T]\times \overline{B(\mathbf{x}_0,R)}\to\mathbb{R}$, $i=1,2,3$, for every $(t,E,H,S)\in [0,24T]\times \overline{B(\mathbf{x}_0,R)}$ by:
\[
f_1(t,E,H,S):=\frac{\eta_0H(P(t,S)-D(t))}{1+\delta_T(T_{\rm ambient}(t)-T_{\rm opt})^2}\left(1-\frac{E}{E_{\rm max}H}\right)-\lambda_0(1+\delta(1-H)^2)E,
            \]
\[
f_2(t,E,H,S):=-\nu\left(\frac{E}{E_{\rm max}H}-50\%\right)^4H,
\]
\[
f_3(t,E,H,S):= \begin{dcases}
(1-S),  & \mbox{if } T_{\rm cell}(t)\geq T_{\text{op}} \\
                   S, & \mbox{otherwise}.
                  \end{dcases}\\
\]
Let us set $\mathbf{g}=(\text{Id},g_2,g_3)$, $\mathbf{x}=(E,H,S)$, and $\mathbf{f}=(f_1,f_2,f_3):[0,24T]\times \overline{B(\mathbf{x}_0,R)}\to\mathbb{R}^3$. Thus, we obtain a concise form for our system:
\begin{equation}\label{eq:concise model}
\begin{aligned}
     &\mathbf{x}_{\mathbf{g}}'(t) = \mathbf{f}(t,\mathbf{x}(t)),\text{ for $\mathbf{g}$-almost every $t\in [0,24T]$,}
     \\
      &\mathbf{x}(0) =\mathbf{x}_0,
      \end{aligned}
\end{equation}
The functions~$T_{\rm ambient}(\cdot)$,~$I(\cdot)$, and $T_{\rm cell}(\cdot)$ are assumed to be continuous, and $D(\cdot)$ is assumed to be measurable thus Lebesgue-Stieltjes measurable. Thus, the function~$\mathbf{f}$ satisfies the assumptions of Theorem~\ref{thm:unique local existence}, and hence there exists $\tau \in (0, 24T]$ such that problem~\eqref{eq:concise model} has a unique solution $\mathbf{x} \in \mathcal{AC}_{\mathbf{g}}([0, \tau], \mathbb{R}^3)$.

\begin{figure}[!h]
\centering
\includegraphics[width=0.9\textwidth]{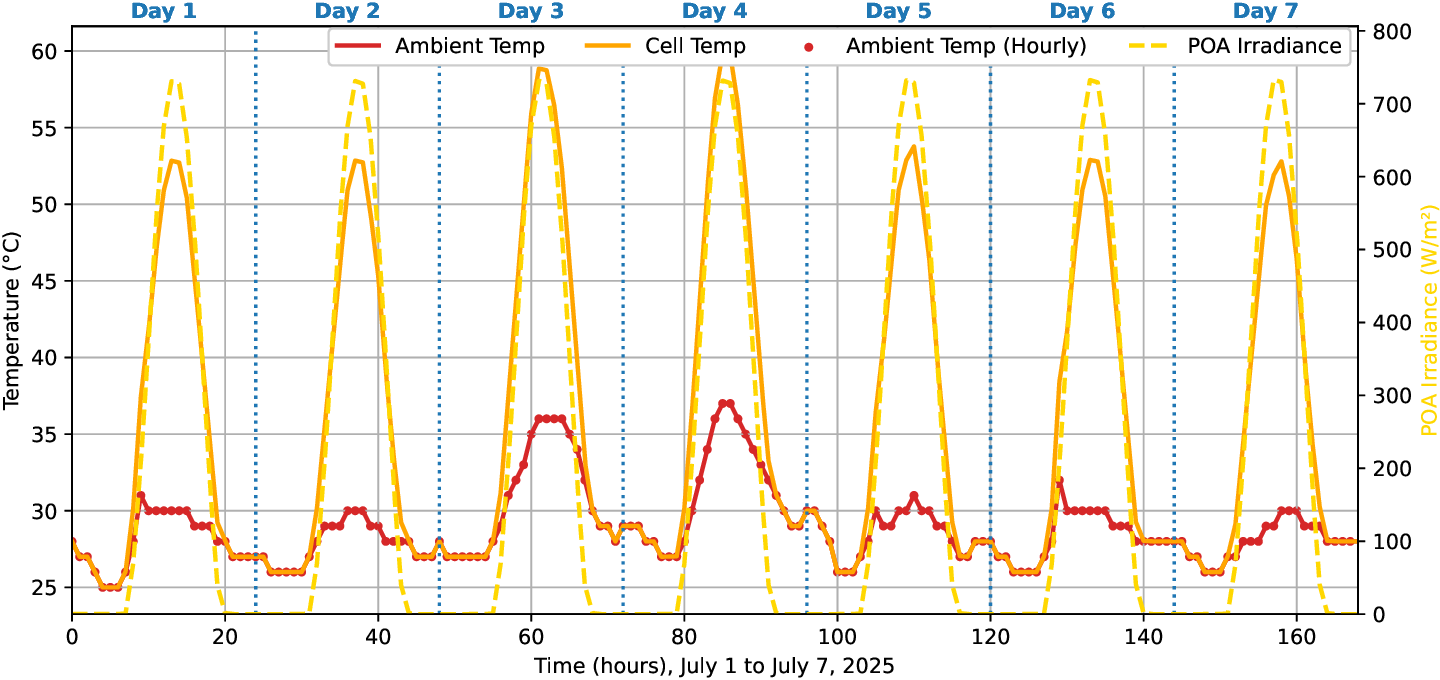}
\caption{Evolution of ambient and cell temperatures and plane-of-array (POA) irradiance over a one-week summer period in Almería, Spain, from 1 July to 7 July 2025. }
\label{fig:Temperatures+Irradiance}
\end{figure}

\paragraph*{Numerical Simulation on a thermal\--stressed interval:}
To illustrate the behavior of the system under thermal stress, we take $T=7$ days and perform a simulation during a summer period—an interval when the ambient temperature remains above the optimal operating temperature $T_{\text{op}}$ of the PV panel, i.e., $T_{\rm cell}(t) > T_{\text{op}}=25^\circ C$ for all $t \in [0,168]$. We set $\beta=1$, which simplifies the derivator~$g_3$ while preserving the essential dynamics:
\[
g_3(t) = \int_{0}^t \mu_1 \cdot \max(0, T_{\rm cell}(s) - T_{\text{op}}) \dif s,
\]
where $T_{\rm cell}(\cdot)$ is shown in Figure~\ref{fig:Temperatures+Irradiance}, using data from Almería for July 1 to July 7, 2025. Ambient temperatures are obtained from local historical weather data~\cite{weatherspark2025}, and the hourly average profiles of direct normal irradiance (DNI) for each month are obtained from the Global Solar Atlas~\cite{globalsolaratlas} and repeated for each day of the 7-day period. Diffuse horizontal irradiance (DHI) and global horizontal irradiance (GHI) are estimated using the Hay\-–Davies decomposition method~\cite{HayDavies1980} to compute the plane\--of\--array (POA) irradiance via \texttt{pvlib} in Python~\cite{pvlib}. Cell temperatures are calculated using a NOCT\--based model, assuming a tilt of $33^\circ$, azimuth of $180^\circ$, NOCT of $45\,^\circ$C, and an albedo of~$0.2$.

In the following, we consider the demand function [W]:
\[
D:t\in [0,168]\to D(t) =
\begin{cases}
120, & 0 \leq t \bmod 24 < 6, \\
180, & 6 \leq t \bmod 24 < 9, \\
130, & 9 \leq t \bmod 24 < 13, \\
180, & 13 \leq t \bmod 24 < 15, \\
130, & 15 \leq t \bmod 24 < 18, \\
200, & 18 \leq t \bmod 24 < 23, \\
120, & 23 \leq t \bmod 24 < 24.
\end{cases}
\]

\begin{figure}[h!]
  \centering
  \includegraphics[width=0.8\textwidth]{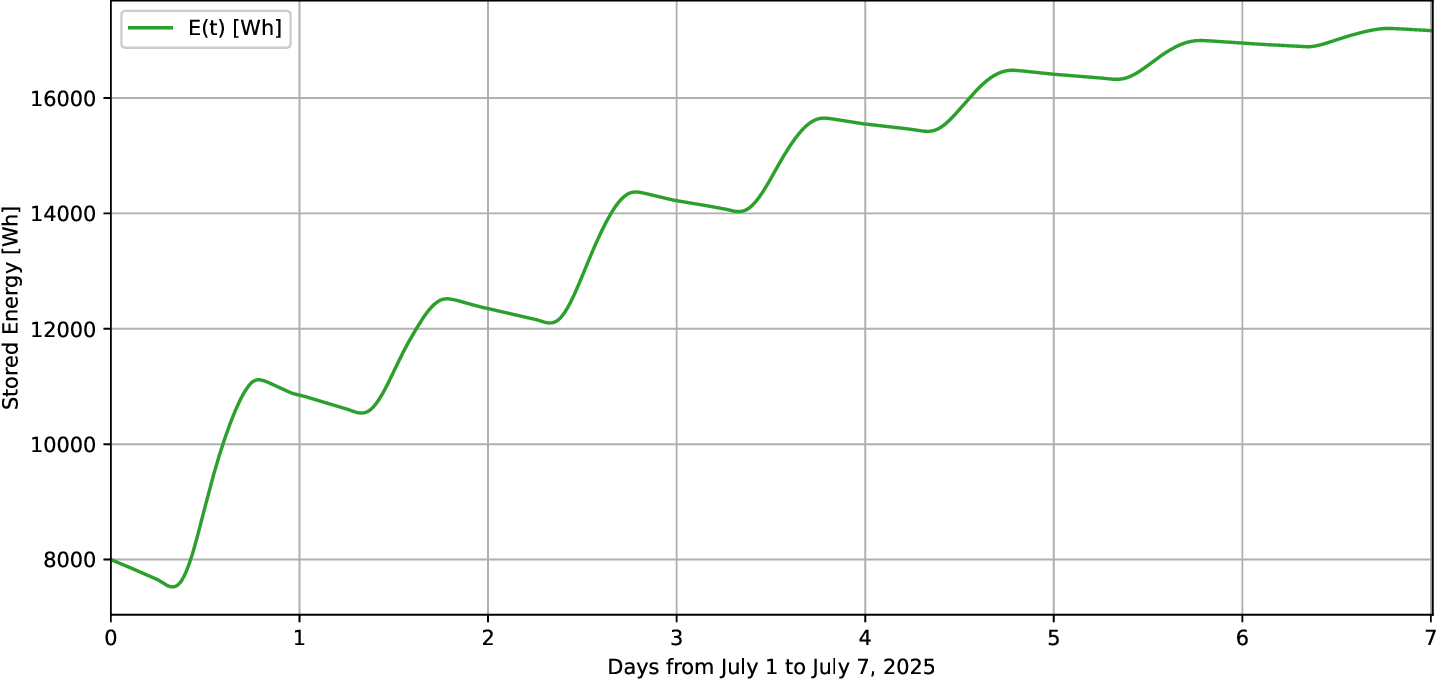}
  \caption{Approximation of battery stored energy with a time-discretization step of 6 minutes over the 7 days.}
  \label{fig:Solution-E}
\end{figure}

\begin{figure}[h!]
  \centering
  \includegraphics[width=0.8\textwidth]{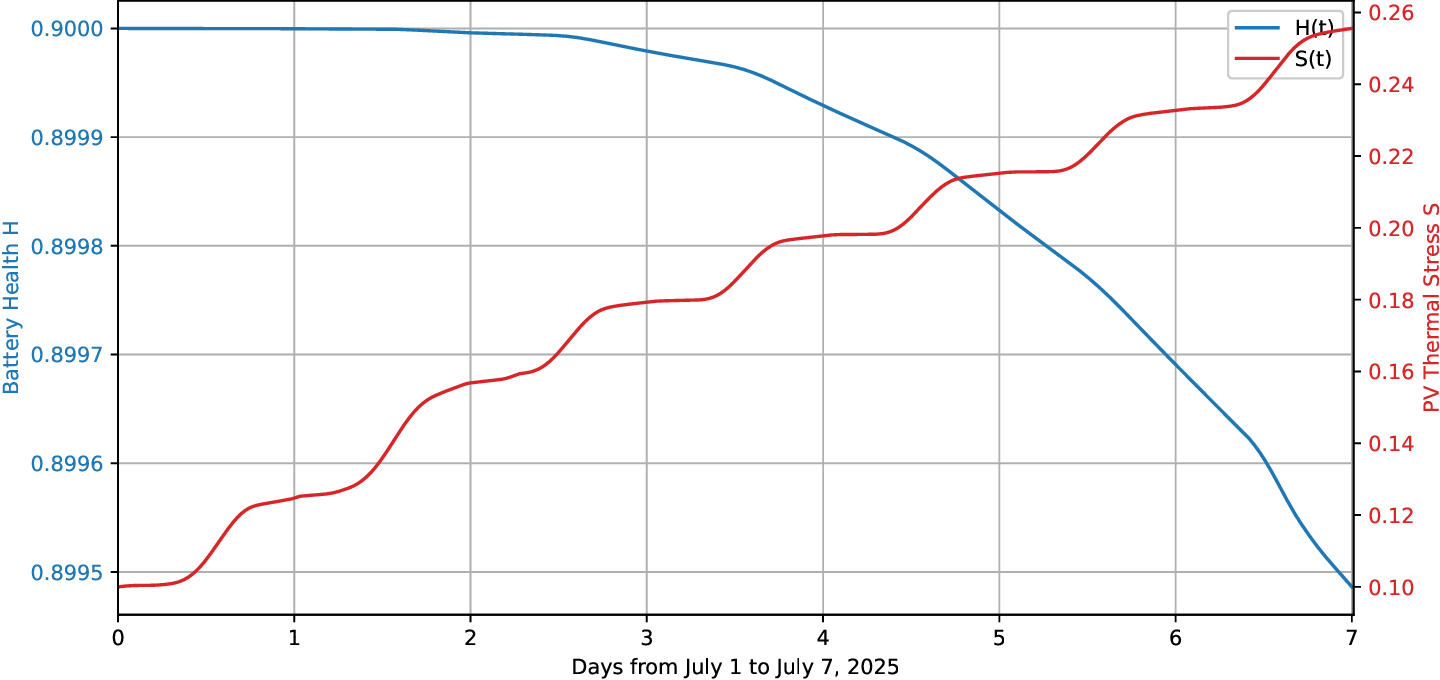}
  \caption{Approximation of battery health and PV thermal stress with a time-discretization step of 6 minutes over the 7 days.}
  \label{fig:Solution-H+S}
\end{figure}

\begin{figure}[ht]
  \centering
  \includegraphics[width=0.8\textwidth]{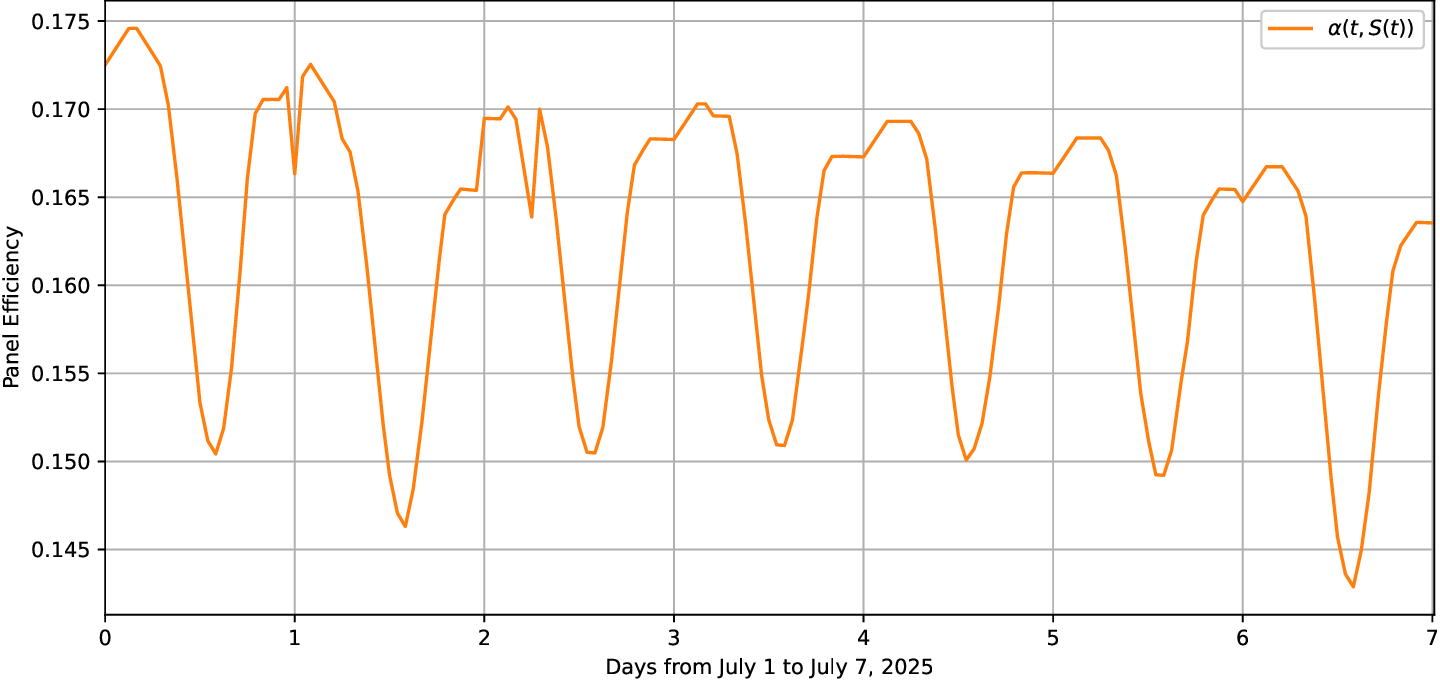}
  \caption{Approximation of the panel efficiency function~$\alpha$.}
  \label{fig:efficiency}
\end{figure}

In Figures~\ref{fig:Solution-E} and~\ref{fig:Solution-H+S}, we provide an approximation of the solution based on the numerical scheme in~\cite{FT}. Parameter values are listed in Tables~\ref{tab:E},~\ref{tab:H}, and~\ref{tab:S}. Figure~\ref{fig:efficiency} shows the panel efficiency function~$\alpha$.

\paragraph*{Comments:}
The numerical simulation result, shown in Figure~\ref{fig:Solution-E} indicates that the stored energy consistently exceeds the instantaneous demand, revealing a surplus that demonstrates the photovoltaic panel’s capacity to meet the load requirements. The stored energy function approaches a short-term equilibrium, reflecting the balance between generation and consumption. As illustrated in Figure~\ref{fig:Solution-H+S}, the battery health function exhibits a gradual decline over time, with the rate of degradation increasing towards the end of the simulated period. Concurrently, the PV thermal stress function steadily rises. This behavior aligns with the system dynamics. The surplus energy leads to higher battery usage, pushing the battery closer to its maximum capacity and accelerating degradation. Meanwhile, persistent cell temperatures above the optimal operating temperature~$T_{\rm op}$ increase thermal stress, limiting recovery and further contributing to the decline in battery health. This behavior is also apparent in Figure~\ref{fig:efficiency}, which illustrates the panel efficiency function~$\alpha$ and shows how the efficiency peak decreases as thermal stress accumulates.


\begin{table}[h!]
\centering
\renewcommand{\arraystretch}{1.5}
\begin{tabular}{|c|c|c|c|c|c|c|}
\hline
$E_0$ [kWh] &$E_{\max}$ [kWh] & $\eta_0$ [--] & $T_{\text{opt}}$ [$^\circ$C] & $\delta_T$ [1/ $^\circ$C$^2$] & $\lambda_0$ [1/h] & $\delta$ [--] \\
\hline
0.8\,$E_{\max}$ &20 & 0.95 & 25 & 0.005 & 0.0001 & 1 \\
\hline
\end{tabular}
\caption{Parameters for the energy storage equation.}\label{tab:E}
\end{table}
\begin{table}[h!]
\centering
\renewcommand{\arraystretch}{1.5}
\begin{tabular}{|c|c|c|c|c|}
\hline
$H_0$ [--] & $\nu$ [1/h] & $T_{\text{thresh}}$ [$^\circ$C] & $\beta_{\text{thermal}}$ [1/$^\circ$C] \\
\hline
0.9 & 0.0002 & 30 & 0.07 \\
\hline
\end{tabular}
\caption{Parameters for the battery health degradation equation.}\label{tab:H}
\end{table}
\begin{table}[h!]
\centering
\renewcommand{\arraystretch}{1.5}
\begin{tabular}{|c|c|c|c|c|c|c|c|}
\hline
$S_0$ [--] & $A$ [m$^2$] & $\alpha_{\text{ref}}$ [--] & $\gamma$ [1/$^\circ$C] & $\rho$ [--] &$\mu_1$ [1/($^\circ C \cdot$h)] & $T_{\text{opt}}$ [$^\circ$C] \\
\hline
0.1  &18 & 0.18 & 0.004 & 0.3 &  0.0001 & 25  \\
\hline
\end{tabular}
\caption{Parameters for PV panel efficiency and thermal stress accumulation.}\label{tab:S}
\end{table}

\section*{Acknowledgements}
Lamiae Maia would like to express her sincere gratitude to Professor Fernando Adrián Fernández Tojo and to the Departamento de Estatística, Análise Matemática e Optimización at the Universidade de Santiago de Compostela, for their warm hospitality during her research stay  at the
 aforementioned department. Lamiae Maia also acknowledges the Galician Center for Mathematical Research and Technology (CITMAga) for funding this stay, during which the present article was finalized.

\section*{Funding}
F. Adrián F. Tojo was supported by grant PID2020-113275GB-I00 funded by Xunta de Galicia, Spain, project ED431C 2023/12; and by MCIN/AEI/10.13039/ 501100011033, Spain, and by “ERDF A way of making Europe” of the “European Union”.

\bibliography{FullBibGD}
\bibliographystyle{spmpsciper}

\end{document}